\documentclass[10pt,a4paper]{amsart}
\usepackage{verbatim}

\usepackage[T1]{fontenc}
\usepackage{graphicx}
\usepackage{enumerate}
\usepackage{amsmath,amsfonts,amssymb}
\usepackage{mathtools}
\usepackage{ulem}
\usepackage{color}
\usepackage{cite}
\definecolor{citation}{rgb}{0.2,0.58,0.2} 
\definecolor{formula}{rgb}{0.1,0.2,0.6}
\definecolor{url}{rgb}{0.3,0,0.5}
\usepackage{pgf,tikz}
\usepackage{mathrsfs}
\usepackage{fancyhdr}
\usepackage{dsfont}

\newcommand{\nocontentsline}[3]{}
\newcommand{\tocless}[2]{\bgroup\let\addcontentsline=\nocontentsline#1{#2}\egroup}

\title[Wolff potentials and vectorial problems with  Orlicz growth]{Wolff potentials and measure data vectorial problems with  Orlicz growth}

\author{Iwona Chlebicka}\address{Iwona Chlebicka \\
Institute of Applied Mathematics and Mechanics, University of Warsaw \\ ul. Banacha 2, 02-097 Warsaw, Poland\\  \texttt{e-mail: i.chlebicka@mimuw.edu.pl}} 

\author{Yeonghun Youn}\address{Yeonghun Youn\\ Department of Mathematics, Yeungnam University\\ Gyeongbuk 38541, Republic of Korea\\ \texttt{e-mail: yeonghunyoun@yu.ac.kr}}

\author{Anna Zatorska-Goldstein}\address{Anna Zatorska-Goldstein \\ Institute of Applied Mathematics and Mechanics, University of Warsaw \\ ul. Banacha 2, 02-097 Warsaw, Poland \\ \texttt{e-mail: azator@mimuw.edu.pl}}

\date{}
\begin{document}

\thanks{\vspace{1mm}}
\maketitle \sloppy

\thispagestyle{empty}

\belowdisplayskip=18pt plus 6pt minus 12pt \abovedisplayskip=18pt
plus 6pt minus 12pt
\parskip 4pt plus 1pt
\parindent 0pt

\newcommand{\ic}[1]{\textcolor{teal}{#1}}

\def\tens#1{\pmb{\mathsf{#1}}}
\newcommand{\barint}{
         \rule[.036in]{.12in}{.009in}\kern-.16in
          \displaystyle\int  } 
          
\newcommand{\dv}{{\rm div}}
\def\aI{\texttt{(a1)}}
\def\aII{\texttt{(a2)}}
\newcommand{\opA}{{\mathcal{ A}}}
\newcommand{\bopA}{{\bar{\opA}}}
\newcommand{\wt}{\widetilde}
\newcommand{\ve}{\varepsilon}
\newcommand{\vp}{\varphi}
\newcommand{\vt}{\vartheta}
\newcommand{\vr}{\varrho}
\newcommand{\pa}{\partial}
\newcommand{\cW}{{\mathcal{W}}}
\newcommand{\supp}{{\rm supp}}

\def\R{{\mathbb{R}}}
\def\N{{\mathbb{N}}}
\def\rp{{[0,\infty)}}
\def\r{{\mathbb{R}}}
\def\n{{\mathbb{N}}}
\def\l{{\mathbf{l}}}
\def\bu{{\bar{u}}}
\def\bg{{\bar{g}}}
\def\bG{{\bar{G}}}
\def\ba{{\bar{a}}}
\def\bv{{\bar{v}}}
\def\wtgamma{{\wt\gamma}}
\def\tew{{\tens{w}}}
\def\teeta{{\tens{\eta}}}
\def\texi{{\tens{\xi}}}
\def\teu{{\tens{u}}}
\def\tebu{{\bar{\tens{u}}}}
\def\teet{{\tens{\eta}}}
\def\teph{{\tens{\phi}}}
\def\tev{{\tens{v}}}
\def\telambda{{\tens{\lambda}}}
\def\calV{{\mathcal{V}}}
\def\tebv{{\bar{\tens{{v}}}}}
\def\tevp{{\tens{\vp}}}
\def\teF{{\tens{F}}}
\def\tef{{\tens{f}}}
\def\teg{{\tens{g}}}
\def\teelvr{{\tens{\ell_\vr}}}
\def\tewtu{{\tens{\wt u}}}
\def\tfA{{{\mathfrak{A}}}}

\def\teDu{{D\teu}}
\def\teDbu{{D\tebu}}
\def\teDet{{D\teet}}
\def\teDph{{D\teph}}
\def\teDv{{D\tev}}
\def\teDw{{D\tew}}
\def\teDbv{{D\tebv}}
\def\teDwtu{{D\tewtu}}
\def\teDvp{D{\tevp}}
\def\teDelvr{{D\teelvr}}

\def\tedv{{\tens{\dv}}}
\def\temu{{\tens{\mu}}}
\def\tea{\tens{a}}

\def\tI{\text{I}}
\def\tII{\text{II}}
\def\tIII{\text{III}}
\def\bmu{{\bar{\mu}}}
\def\rn{{\mathbb{R}^{n}}}
\def\Rm{{\mathbb{R}^{m}}}
\def\Rn{{\mathbb{R}^{n}}}
\def\id{{\mathsf{Id}}}
\def\P{{\mathsf{P}}}
\def\Pj{{\mathsf{P_j}}}
\def\rnm{{\mathbb{R}^{n\times m}}}
\def\rN{{\mathbb{R}^{N}}} 
\def\Mb{{\mathcal{M}(\Omega,\Rm)}} 

\newtheorem{coro}{\bf Corollary}[section]
\newtheorem{theo}[coro]{\bf Theorem} 
\newtheorem{lem}[coro]{\bf Lemma}
\newtheorem{rem}[coro]{\bf Remark} 
\newtheorem{defi}[coro]{\bf Definition} 
\newtheorem{ex}[coro]{\bf Example} 
\newtheorem{fact}[coro]{\bf Fact} 
\newtheorem{prop}[coro]{\bf Proposition}

\newcommand{\data}{\textit{\texttt{data}}}


\parindent 1em

\begin{abstract}
We study  solutions to measure data elliptic systems with Uhlenbeck-type structure   
that involve operator of divergence form, depending continuously on the spacial variable, and exposing doubling Orlicz growth with respect to the second variable. 
Pointwise estimates for the solutions that we provide are expressed in terms of a nonlinear potential of generalized Wolff type.  
Not only we retrieve the recent sharp results proven for $p$-Laplace systems, but additionally our study covers the natural scope of operators with similar structure and natural class of Orlicz growth.
\end{abstract}

\setcounter{tocdepth}{1}
\tableofcontents

\section{Introduction}

A broad and profuse branch of the theory of nonlinear partial differential systems stems from seminal ideas by Ural'tseva~\cite{Ural} and Uhlenbeck~\cite{Uhl}. There is a solid stream of recent studies on regularity of solutions to general growth elliptic systems and related minimizers of vectorial functionals ~\cite{bemi,CiMa-ARMA2014,DiEt,DiMa,DF1,DFMin,DSV1,DSV2,GPT,lieb,Marc2006,MarcPapi,str}. Our aim is to contribute to the field by providing pointwise estimates for very weak solutions to measure data problems in the terms of potentials of relevant Orlicz growth. As a consequence, we infer sharp description of fine properties of solutions being the exact analogues of the ones available in the classical linear potential theory or in the  case of $p$-Laplace systems~\cite{KuMi2018}.

 Let us stress that {weak solutions do} not have to exist for arbitrary measure datum. Thus we employ a notion of very weak solutions obtained by an approximation studied since~\cite{BG}. Despite they 
 can be unbounded, they can be controlled by a certain  potential.   There are known deep classical results for scalar problems~\cite{KiMa92,KiMa94} settling  the nonlinear potential theory for solutions to $-\Delta_p u=-\dv (|Du|^{p-2}Du)=\mu$, $1<p<\infty$, followed by their variable exponent version~\cite{LuMaMa} and recently by the Orlicz one~\cite{CGZG-Wolff,Maly-Orlicz}, as well as a counterpart proven for systems involving $p$-Laplace operator~\cite{CiSch,KuMi2018}. We study the nonstandard growth version of pointwise estimates involving suitably generalized potential of the Wolff type and infer their regularity consequences. 

In fact, we investigate very weak solutions ${\teu}:\Omega\to\Rm$  to measure data elliptic systems involving nonlinear operators 
 \begin{equation}
\label{intro:eq:main} -{ \tedv}\left( {a(x)}\frac{g(|\teDu|)}{|\teDu|}  \teDu \right)=\temu\quad\text{in}\quad \Omega,
\end{equation}
where $a\in C(\Omega)$ is bounded and separated from zero, $g=G'$, $\temu\in\mathcal{M}(\Omega,\Rm)$ is a bounded measure, whereas $\tedv$ stands for the $\Rm$-valued divergence
operator. Here we admit  $G$ to be a Young function satisfying both $\Delta_2$- and $\nabla_2$-conditions, which follows from the condition 
\begin{equation}\label{iG-sG} 
2 \leq i_G=\inf_{t>0}\frac{tg(t)}{G(t)}\leq \sup_{t>0}\frac{tg(t)}{G(t)}=s_G<\infty.\end{equation}  
See {\it Assumption }{\bf (A-vect)} in Section~\ref{sec:formulation-n-reg-cons} for all details.   We stress that we impose the typical assumption of quasi-diagonal structure of  $\opA$ naturally covering the case of (possibly weighted) $p$-Laplacian when $G_p(s)=s^{p},$ for every $p\geq 2$, together with operators governed by the Zygmund-type functions $G_{p,\alpha}(s)=s^p\log^\alpha(1+s)$, $p\geq 2,\,\alpha\in \R$, as well as their multiplications and compositions with various parameters. In turn, we generalize the corresponding results of~\cite{CiSch,KuMi2018} to embrace also $p$-Laplace systems with continuous coefficients, and -- on the other hand -- cover the natural scope of operators with similar structure and Orlicz growth. Let us stress that the operator we consider is {\em not} assumed to enjoy homogeneity of a~form $\opA(x,k\xi)=|k|^{p-2}k\opA(x,\xi)$. Consequently, our class of solutions is {\em not} invariant with respect to scalar multiplication.

Obtaining sharp regularity results for solutions to nonlinear systems is particularly challenging. In the scalar case one can infer continuity or H\"older continuity of the solution to $-\dv \opA(x,Du)=\mu$ when the dependence of the operator on the spacial variable is merely bounded and measurable and the growth of $\opA$ with respect to the second variable is governed by an arbitrary doubling Young function (cf. \cite{CGZG-Wolff}). The same is not possible for systems even with far less complicated growth and null datum, cf.~\cite{DeG,JMVS,SY} and~\cite[Section~3]{Min-Dark}.  To justify why we restrict our attention to operators having a specific form as in~\eqref{intro:eq:main}, let us point out that the typical assumption of a so-called Uhlenbeck structure is imposed in order to control energy of solutions. The continuity of the coefficients is a minimal assumption to get continuity of the solution in the view of counterexample of~\cite{DeG}. 

The studies on the potential theory to measure data problems dates back to \cite{HaMa2,Ma}. We refer to~\cite{KuMi2014} for an overview of the nonlinear potential theory, to \cite{HedWol,KiMa92,KiMa94,LiMa} for cornerstones of the field, and~\cite{adams-hedberg,hekima} for well-present background of the $p$-growth case. In the scalar case, in their seminal works Kilpel\"ainen and Mal\'y in~\cite{KiMa92,KiMa94} provided optimal Wolff potential estimates for $p$-superharmonic functions $u$ generating a nonnegative measure $\mu$ from above and below
\begin{equation}
    \label{wolff-p-est}
\tfrac 1c\cW_p^{\mu}(x_0,R)\leq u(x_0)\leq c\Big(\inf_{B_R(x_0)}u+\cW_p^{\mu}(x_0,R)\Big)\quad\text{for some }\ c=c(n,p)
\end{equation}with {the} so-called Wolff potential
\[ \cW_p^{\mu}(x_0,R)=\int_0^R \left(r^{p-n}{\mu(B_r(x_0))}\right)^{\frac{1}{p-1}}\,\frac{dr}{r}= \int_0^R \left(\frac{\mu(B_r(x_0))}{r^{n-1}}\right)^{\frac{1}{p-1}}\,dr,\]
see also~\cite{KoKu,tru-wa}. In the linear case ($p=2$) the estimates of~\eqref{wolff-p-est} become the classical Riesz potential bounds. The precise Orlicz counterpart of this result with nonnegative measure $\mu$ is proven with the nonstandard growth potential \begin{equation*}
\cW^{\mu}_G(x_0,R)=\int_0^R g^{-1}\left(\frac{\mu(B_r(x_0))}{r^{n-1}}\right)\,dr,
\end{equation*}
see \cite{CGZG-Wolff,Maly-Orlicz}. To our best knowledge, the only reference one can find on the related results for systems are~\cite{CiSch,KuMi2018} that involves problems with the $p$-Laplace operator. The estimate related to~\eqref{wolff-p-est} provided therein establishes the upper bound only. Note however that no lower bound can be available in the vectorial case. Indeed, it origins in the lack of the possibility of proving maximum principle. Our method of proof relies on the ideas of \cite{KuMi2018}. We employ a properly adapted Orlicz version of $\opA$-harmonic approximation relevant for measure data problems (Theorem~\ref{theo:Ah-approx} in Section~\ref{sec:Ah-approx}) and careful estimates on concentric balls.

Let us also mention that a lot of  attention is attracted by potential estimates on gradients of solutions, generalized harmonic approximation, and their application in the theory of partial regularity \cite{Baroni-Riesz,BaHa,BCDKS,Byun4,Byun5,Byun6,DiLeStVe,DSV3,DuGr,DuMi2010-2,KuMi2016p,KuMi2014} which is an open path from now on.
 
The paper is organized as follows. Our assumptions, main results and their regularity consequences are presented in Section~\ref{sec:formulation-n-reg-cons}. Section~\ref{sec:prelim} is devoted to notation and information on the setting. In particular  see Section~\ref{ssec:sols} for the precise definition of the notion of very weak solutions we employ. Section~\ref{sec:Ah-approx} provides the most important tool of paper -- a measure data $\opA$-harmonic approximation. Section~\ref{sec:mainproof} contains the proofs of comparison estimates, the sufficient condition for $\teu$ to be in VMO of Proposition~\ref{prop:vmo}, the potential estimates of Theorems~\ref{theo:pointwise}, the continuity criterion of Theorem~\ref{theo:continuity} and the H\"older continuity criterion of Theorem~\ref{theo:H-cont}.

\section{Main result and its consequences }\label{sec:formulation-n-reg-cons}

\subsection{The statement of the problem}
Let us present an essential notation and details of the measure data problem we study.

\noindent{\underline{\it Essential notation}.} By {`$\cdot$'} we denote the {scalar product of two vectors, i.e.  for ${\texi}=(\xi_1,\dots,\xi_m)\in \Rm$ and 
	$\teeta= (\eta_1,\dots,\eta_m)\in \Rm$ we have ${\texi}\cdot{\teeta} = \sum_{i=1}^m \xi_i \eta_i$};  by {`$:$'} -- {the Frobenius product of the second-order tensors, i.e. for ${\xi}=[\xi_{j}^\alpha]_{j=1,\dots,n,\, \alpha=1,\dots,m}$ and $\eta=[\eta_{j}^\alpha]_{j=1,\dots,n,\, \alpha=1,\dots,m}$ we have
	\[{\xi}: {\eta} =\sum_{\alpha=1}^m \sum_{j=1}^n \xi_{j}^\alpha \eta_{j}^\alpha.\]}
	By `{$\otimes$}' we denote {the tensor product of two vectors, i.e for ${{\texi}}=(\xi_1,\dots,\xi_k)\in {\R^k}$  and 
	${{\teeta}}= (\eta_1,\dots,\eta_\ell)\in \R^\ell$, we have $\texi\otimes\teeta:=[\xi_i\eta_j]_{i=1,\dots,k,\,j=1,\dots,\ell},$ that is \[{\texi}\otimes{\teeta}:=\begin{pmatrix}
  \xi_{1}\eta_{1} & \xi_{1}\eta_{2} & \cdots & \xi_{1}\eta_{\ell} \\
  \xi_{2}\eta_{1} & \xi_{2}\eta_{2} & \cdots & \xi_{2}\eta_{\ell} \\
  \vdots  & \vdots  &  & \vdots  \\
  \xi_{k}\eta_{1} & \xi_{k}\eta_{2} & \cdots & \xi_{k}\eta_{\ell} 
 \end{pmatrix}\in \R^{k\times \ell}.\]}

\medskip

\noindent{\underline{\it Assumption {\bf (A-vect)}}.}  Given a bounded, open, {Lipschitz} set $\Omega\subset\rn$, $n\geq 2$, we investigate solutions ${\teu}:\Omega\to\Rm$ to the problem \begin{equation}
    \label{eq:mu}\begin{cases}-\tedv \opA(x,\teDu)=\temu\quad\text{in }\ \Omega,\\
    {\teu}=0\quad\text{on }\ \partial\Omega\end{cases}
\end{equation}
 with a datum $\temu$ being a vector-valued bounded Radon measure and a function $\opA:\Omega\times\rnm\to\rnm$ is a weighted operator of Orlicz growth expressed by the means of $g(t):=G'(t)$, where an $N$-function $G \in C^2((0,\infty))\cap C(\rp)$ satisfies $i_G\geq 2$ with $i_G$ given by~\eqref{iG-sG}. {Let $g\in\Delta_2\cap\nabla_2$.}  Namely, $\opA$ is assumed to admit a~form\begin{equation}
    \label{opA:def}\opA(x,\xi)=a(x)\frac{g(|\xi|)}{|\xi|}\,\xi,\end{equation} where $a:\Omega\to[c_a,C_a],$ $0<c_a<C_a$ is a continuous function  with a modulus of continuity $\omega_a$. 
 We define a potential
  \begin{equation}
\label{Wolff-potential}
\mathcal{W}^{\mu}_G(x_0,R)=\int_0^R g^{-1}\left(\frac{|\temu|(B_r(x_0))}{r^{n-1}}\right)\,dr.
\end{equation} For the case of referring to the dependence of some quantities on the parameters of the problem, we collect them as \[\data=\data(i_G,s_G,c_a,C_a,\omega_a,n,m).\]Having~\eqref{opA:def}, one can infer the strong monotonicity of the vector field $\opA$ of a form given by Lemma~\ref{lem:DiEt-mon}. 

\medskip

Let us define the notion of very weak solutions we employ. 

\noindent{\underline{\it Solutions Obtained as a Limit of Approximation}.}  A map $\teu\in W^{1,g}(\Omega,\Rm)$  is called a SOLA to~\eqref{eq:mu} under the regime of {\rm Assumption {\bf (A-vect)}}, if there exists a sequence $(\teu_{h})\subset W^{1,G }(\Omega,\Rm)$ of local energy solutions to the systems
\[-\tedv\opA(x,D\teu_{h})={\temu_h}\]
such that $\teu_{h}\to \teu$ locally in $W^{1,g}(\Omega,\Rm)$ and $(\temu_h)\subset C^\infty (\Omega,\Rm)$ is a sequence
of smooth maps that converges to $\temu$ weakly in the sense of measures and satisfies
\begin{equation}
    \label{conv-of-meas}
\limsup_h |\temu_h |(B) \leq |\temu|(B)
\end{equation}
for every ball $B\subset\Omega$. 

Observe that the above approximation property immediately implies that a SOLA $\teu$ is a distributional solution to~\eqref{eq:mu}, that is,
\[\int_\Omega\opA(x,\teDu):\teDvp\,dx=\int\tevp \,d\temu\qquad\text{for every }\tevp\in  C^\infty (\Omega,\Rm).\]

\subsection{Main results}\label{ssec:main-results} Our main accomplishment reads as follows.

  \begin{theo}[Pointwise Wolff potential estimates]\label{theo:pointwise}
    Suppose $\teu:\Omega\to\Rm$ is a SOLA to~\eqref{eq:mu} with $\opA$ satisfying {\rm Assumption {\bf (A-vect)}} and $\temu\in\Mb$. Let $B_r(x_0)\Subset\Omega$ with $r<R_0$ for some $R_0=R_0(\data)$.  If $\cW^\temu_G(x_0,r)$ is finite, then $x_0$ is a Lebesgue's point of $\teu$ and\begin{equation}
        \label{Wolff-osc-est}
        |\teu(x_0)-(\teu)_{B_r(x_0)}|\leq C_\cW\left(\cW^\temu_G(x_0,r)+\barint_{B_r(x_0)}|\teu-(\teu)_{B_r(x_0)}|\,dx\right)
    \end{equation}
holds for $C_\cW>0$ depending only on $\data$. In particular, we have the following pointwise estimate 
\begin{flalign}
\label{eq:u-est}  |\teu(x_0 )|\leq C_\cW\left(\mathcal{W}^{\temu}_G (x_0 ,r)+\barint_{B_r(x_0)} |\teu(x)|dx\right).
\end{flalign}
\end{theo}
 
 \subsection{Local behaviour of very weak solutions} \label{ssec:loc-beh} Potential estimates are known to be an efficient tool to bring precise information on the local behaviour of solutions. We refer to~\cite{KuMi2014} for clearly presented overview of consequences of estimates like~\eqref{eq:u-est} in studies on $p$-superharmonic functions together with a bunch of related references and to~\cite{CGZG-Wolff} for similar results for $\mathcal{A}$-harmonic functions where the operator $\mathcal{A}$ is exhibiting Orlicz type of growth. Notice however we referred to the scalar case. In the vectorial one, the only investigations on the potential estimates to solutions to measure data problem we are aware of are available in~\cite{CiSch,KuMi2018} for $p$-Laplace systems. 
 Let us present the regularity consequences of  Wolff potential estimates to $p$-Laplace system with continuous  coefficients and the natural scope of operators with similar structure and Orlicz growth.
 
 \medskip
 
We start with finding a density condition around a point $x_0$ implying that a solution has vanishing mean oscillations at $x_0$. The proposition below does not follow from Theorem \ref{theo:pointwise}, for the proof see Section~\ref{sec:mainproof}.
 
   \begin{prop}[VMO criterion] \label{prop:vmo}  Suppose $\opA$ satisfies {\rm Assumption {\bf (A-vect)}} and $\temu\in\Mb$.      Let $\teu$ be a SOLA to~\eqref{eq:mu} and let $B_r(x_0)\Subset\Omega.$ If\begin{equation}
         \label{mu-shrinks}\lim_{\vr\to 0}\vr\, g^{-1}\left(\frac{|\temu|(B_\vr(x_0))}{\vr^{n-1}}\right)=0,
     \end{equation}
     then $\teu$ has vanishing mean oscillations at $x_0$, i.e.
     \begin{equation}
         \label{u-VMO-x0}\lim_{\vr\to 0}\,\barint_{B_\vr(x_0)}|\teu-(\teu)_{B_\vr(x_0)}|\,dx=0.
     \end{equation}
 \end{prop} 
 
An application of Theorem~\ref{theo:pointwise} is the following continuity criterion proven also in Section~\ref{sec:mainproof}.
 
 \begin{theo}[Continuity criterion]
     \label{theo:continuity} Suppose $\teu$ be a SOLA to~\eqref{eq:mu} under the regime of {\rm Assumption {\bf (A-vect)}} and $B_r(x_0)\Subset\Omega.$ If\begin{equation}
         \label{Wolff-shrinks}
         \lim_{\vr\to 0}\sup_{x\in B_r(x_0)}\cW^\temu_G(x,\vr)=0,
     \end{equation}
then $\teu$ is continuous in $B_r(x_0).$
 \end{theo}

If $\temu=0$, then  $\cW^\temu_G(x,\vr)=0$, thus trivially we have the following consequence.
\begin{coro}\label{coro:ah-cont}
Under {\rm Assumption {\bf (A-vect)}} if $\teu$ is an $\opA$-harmonic map in $\Omega,$ then $\teu$ is continuous in every $\Omega'\Subset\Omega.$ 
\end{coro}
 
Condition~\eqref{Wolff-shrinks} holds true provided the datum belongs to a Lorentz-type space. In order to define it we recall some definitions. We denote by $f^\star$  the decreasing rearrangement of a measurable function $f:\Omega\to\R$ by
$$
f^\star(t) = \sup \{ s \geq 0 \colon |\{x\in \R^n:f(x)>s\}| > t  \},
$$
the maximal rearrangement by
$$
f^{\star \star}(t) =\frac 1t \int_0^t f^\star(s) \, ds\quad\text{and}\quad 
f^{\star \star}(0)=
f^{\star}(0).
$$
Following \cite{stein-w} by Lorentz space $L(\alpha, \beta)(\Omega)$ for $\alpha,\beta>0$ we mean the class of measurable functions such that  
$$ 
 \int_0^\infty \left( t^{1/ \alpha} f^{\star \star}(t)\right)^\beta\,\frac{dt}{t} 
 <\infty.
$$ 
{The following fact is proven in Appendix.}
\begin{lem}\label{lem:Wolff-est} Suppose  $\temu=\teF:\Rn\to\Rm$ is a locally integrable map vanishing outside $\Omega$, then there exists a constant $c=c(n,i_G,s_G)>0$, such that\[\cW_G^{\teF}(x,R)\leq c\int_0^{|B_R|}  t^\frac{1}{n} g^{-1}\left(t^\frac{1}{n} |\teF|^{\star\star}(x)\right)\,\frac{dt}{t}=: \mathcal{I}_R.\]
\end{lem}
Note that if $\mathcal{I}_R<\infty$, the datum $\teF$ is in the dual space to $W^{1,G}(B_R)$ and, consequently, we deal with a weak solution, see~\cite{ACCZG}. Let us present a corollary of Theorem~\ref{theo:continuity} holding due to Lemma~\ref{lem:Wolff-est}. 
 \begin{coro}
 \label{coro:cont-2} If $\teu$ is a weak solution to $-\tedv\opA(x,\teDu)=\teF$ with $\opA$ satisfying {\rm Assumption {\bf (A-vect)}} and $\teF:\Omega\to\Rm$ such that \begin{equation}
    \label{cond-Lor}
 \mathcal{I}:=\int_0^{|\Omega|} t^\frac{1}{n} g^{-1}\left(t^\frac{1}{n} |\teF|^{\star\star}(x)\right)\,\frac{dt}{t}  <\infty\end{equation}
for $\Omega_0\Subset\Omega$, then $\teu\in C(\Omega_0,\Rm)$ {and $\|\teu\|_{L^\infty(\Omega_0,\Rm)}\leq c(\data) \mathcal{I}.$ This bound is optimal and attained by a radial solution on a ball, see \cite{ACCZG}. Moreover, as special cases we get that $\teu$ is continuous under the regularity restrictions on $|\teF|=f$ of \cite[Example~1 (A) and Example~2 (A)]{ACCZG}, still within our regime requiring $g\in\Delta_2\cap\nabla_2$ and $i_G\geq 2$.}\end{coro}
The above corollary 
results in the following extension of \cite[Theorem 10.6]{KuMi2018} to the weighted case.

\begin{rem}\label{rem:p-Lorentz}\rm If $\teu$ is a weak solution to $-{ \tedv}\left(a(x)|\teDu|^{p-2} \teDu \right)=\teF$ for $p\geq 2$  and $0<a\in C(\Omega)$ is separated from zero and $|\teF|$ belongs locally to the  Lorentz space $L(\tfrac{n}{p},\tfrac{1}{p-1})(\Omega)$ then $\teu$ is continuous in $\Omega$. \end{rem} As another application of Corollary~\ref{coro:cont-2} let us present its consequences for the Zygmund case.
\begin{rem} \rm Suppose that $2\leq p<n,$ {$\alpha\geq 0$}, $0<a\in C(\Omega)$  separated from zero is bounded, and $\teu$ is a weak solution to \begin{equation}
    \label{eq-log-Fu} -{ \tedv}\left(a(x)|\teDu|^{p-2}\log^\alpha({\rm e}+|\teDu|)  \,\teDu \right)=\teF.\end{equation} Observe that in this case $g^{-1}(\lambda)\approx \lambda^\frac{1}{p-1}\log^{-\frac{\alpha}{p-1}}({\rm e}+\lambda)$. If $\teF$ satisfies~\eqref{cond-Lor}, then $\teu$ is continuous in $\Omega$. 
\end{rem}

  Wolff potential estimates can be used to find a relevant conditions on a measure $\temu$ to infer H\"older continuity of solutions. One of the natural ones is expressed in the Orlicz modification of the Morrey-type scale. \begin{theo}[H\"older continuity criterion] \label{theo:H-cont}  Suppose $\teu$ be a SOLA to~\eqref{eq:mu} under the regime of {\rm Assumption {\bf (A-vect)}}. Assume further that for $\temu$ there exist positive constants $c=c(\data)>0$ and $\theta\in (0,1)$ such that  \begin{equation}
    \label{mu-control}
|\temu|(B_r(x))\leq c r^{n-1}g(r^{\theta-1})\approx r^{n-\theta}G(r^{\theta-1})
\end{equation}
for each $B_r(x)\Subset\Omega$ with sufficiently small radius.   Then   $\teu$ is locally H\"older continuous in $\Omega$.  
\end{theo}   For $p$-growth problems condition~\eqref{mu-control} reads as $|\temu|(B_r(x))\leq c  r^{n-p+\theta(p-1)}$  well known since~\cite{Car,KiMa92,kizo,RaZi}. Moreover, in the scalar case~\eqref{mu-control} is proven in~\cite{CGZG-Wolff} to be equivalent to H\"older continuity of solutions, while in~\cite{ChKa} to characterize removable sets for H\"older continuous solutions.  In the vectorial case we cannot get equivalence, because by Theorem~\ref{theo:pointwise} we are equipped with one-sided estimate only. Specializing Theorem~\ref{theo:H-cont}, we have the following results.
\begin{rem} \rm Suppose  $p\geq 2$, positive $a\in C(\Omega)$ is separated from zero, and $\teu$ is a SOLA to \[ -{ \tedv}\left(a(x)|\teDu|^{p-2} \teDu \right)=\temu\]  with  $|\temu|(B_r(x))\leq c r^{n-p+\theta(p-1)}$ 
for some $c>0,$ $\theta\in(0,1)$ and all sufficiently small $r>0.$ Then $\teu$ is locally  H\"older continuous in $\Omega$.
\end{rem}
\begin{rem} \rm Suppose  $p\geq 2,$ $\alpha\in\R$, positive $a\in C(\Omega)$ is separated from zero, and $\teu$ is a SOLA to  
\begin{equation}
    \label{eq-log-mu}- \tedv\left(a(x)|\teDu|^{p-2}\log^\alpha({\rm e}+|\teDu|)  \teDu \right)=\temu
\end{equation}   with  $|\temu|(B_r(x))\leq c r^{n-p+\theta(p-1)}\log^\alpha({\rm e}+r^{\theta-1}) $ for some $c>0,$ $\theta\in(0,1)$ and all sufficiently small $r>0.$ Then $\teu$ is locally  H\"older continuous in $\Omega$.
\end{rem}

The sufficient condition for~\eqref{mu-control} and, in turn, for the H\"older continuity of the solution is to assume that $|\temu|=|\teF|$ belongs to a~relevant Marcinkiewicz-type space. Following~\cite{oneil} for a continuous increasing function $\psi:(0,|\Omega|)\to (0,\infty)$ we say that $f\in L(\psi,\infty)(\Omega)$ if the maximal rearrangement $f^{\star\star}$ of $f$ satisfies
\[\sup_{s\in(0,|\Omega|)}\frac{f^{\star\star}(s)}{\psi^{-1}(1/s)}<\infty.\]
There holds the following consequence of Theorem~\ref{theo:H-cont}.

\begin{coro}\label{coro:O-Marc-dens}  Suppose $\teu$ is a SOLA to $-\tedv\opA(x,\teDu)=\teF$  under the regime of {\rm Assumption {\bf (A-vect)}}   and $|\teF|$ belongs locally to the Marcinkiewicz-type space $L(\psi,\infty)(\Omega)$ with $\psi^{-1}(1/\lambda)=\lambda^{-\frac{1}{n}}g\big(\lambda^{\frac{\theta-1}{n}}\big)$ for some $\theta\in(0,1)$, 
then $\teu$ is locally H\"older continuous. \end{coro} 
For justification that indeed under the assumptions of Corollary~\ref{coro:O-Marc-dens} the condition~\eqref{mu-control} is satisfied see the calculations provided for the scalar case \cite[Section~2]{CGZG-Wolff}. The above fact has  the best possible consequence in the $p$-Laplace case.
\begin{rem}\rm \label{rem:p-Marc-dens} If $p\geq 2$, positive $a\in C(\Omega)$ is separated from zero, $\teu$ is a SOLA to $-{ \tedv}\left(a(x)|\teDu|^{p-2} \teDu \right)=\teF$ and $|\teF|$ belongs locally to the Marcinkiewicz space $L(\frac{n}{p+\theta(p-1)},\infty)(\Omega)$  for some $\theta\in(0,1)$, i.e.  $\sup_{\lambda>0}\left(\lambda^\frac{n}{p+\theta(p-1)}\big|\{x\in\Omega_0:\,|\teF(x)|>\lambda\}|\right)<\infty$ for any $\Omega_0\Subset\Omega$, then $\teu$ is locally H\"older continuous in $\Omega$.
\end{rem}
 
\begin{rem}\label{rem:G-Marc-dens} \rm When $G(t)\approx t^p\log^\alpha({\rm e}+t),$  $p\geq 2,$ $\alpha\in\R$, $\teu$ is a SOLA to \eqref{eq-log-Fu} with $\teF$ such that \[\sup_{\lambda>0}\left(\lambda^\frac{n}{p+\theta(p-1)}\log^{-\frac{\alpha(1-\theta)}{p+\theta(p-1)}}({\rm e}+\lambda^\frac{n}{1-\theta})\big|\{x\in\Omega_0:\,|\teF(x)|>\lambda\}|\right)<\infty,\] for any $\Omega_0\Subset\Omega$, then $\teu$ is locally H\"older continuous in $\Omega$.
\end{rem}

\section{Preliminaries}\label{sec:prelim}

\subsection{Notation}
We shall adopt the customary convention of denoting by $c$ a constant that may vary from line to line. To skip rewriting a constant, we use $\lesssim$. By $a\approx b$, we mean $a\lesssim b$ and $b\lesssim a$. To stress the dependence of the intrinsic constants on the parameters of the problem, we write $\lesssim_\data$ or $\approx_\data$. By $B_R$ we denote a ball skipping prescribing its center, when {it} is not important. By $cB_R=B_{cR}$ we mean a ball with the same center as $B_R$, but with rescaled radius $cR$. We make use of the truncation operator,  $T_k:\Rm\to\Rm$, defined as follows 
\begin{equation}\label{Tk}T_k(\texi):=\min\left\{1,\frac{k}{|\texi|}\right\}\texi.
\end{equation}
Then, of course, $D T_k:\Rm\to\R^{m\times m}$ is given by\begin{equation}
    \label{DTk}DT_k(\texi)=\begin{cases}
    \id&\text{if }\ |\texi|\leq k,\\
    \frac{k}{|\texi|}\left(\id-\frac{\texi \otimes \texi}{|\texi|^2}\right)&\text{if }\ |\texi|> k.
    \end{cases}
\end{equation} 

For a~measurable set $U\subset \rn$  with finite and positive $n$-dimensional Lebesgue measure $|U|>0$ and $\tef\colon U\to \mathbb{R}^{k}$, $k\ge 1$ being a measurable map, we define
\begin{flalign*}
(\tef)_U=\barint_{U}\tef(x) \, dx =\frac{1}{|U|}\int_{U}\tef(x) \,dx.
\end{flalign*}
By $C^{0,\gamma}(U)$, $\gamma \in (0,1]$, we mean the family of H\"older continuous functions, i.e. measurable functions $f\colon U\to \mathbb{R}^k$ for which
\begin{flalign*}
[\tef]_{0,\gamma}:=\sup_{\substack{x,y\in U,\\x\not =y}}\frac{|\tef(x)-\tef(y)|}{|{x-y}|^{\gamma}}<\infty.
\end{flalign*} 
We {describe the ellipticity of a vector field $\opA$} using a function $\calV:\rnm\to\rnm$
\begin{equation}
    \label{V-def}
\calV({\xi})= \left(\frac{g(|{\xi}|)}{|{\xi}|}\right)^{1/2}\xi . 
\end{equation}

\subsection{Basic definitions} References for this section {are~\cite{rao-ren,KR}}.

We say that a function $G: [0, \infty) \to [0, \infty]$ {is a Young function} if it is convex, vanishes  at $0$, and is neither identically equal to 0, nor to infinity. A Young function $G$ which is finite-valued,
vanishes only at $0$ and satisfies the additional growth conditions
\begin{equation*}\lim _{t \to
0}\frac{G (t)}{t}=0 \qquad \hbox{and} \qquad \lim _{t \to \infty
}\frac{G (t)}{t}=\infty 
\end{equation*} 
is called an $N$-function.  The  complementary~function $\wt{G}$  (called also the Young conjugate, or the Legendre transform) to a nondecreasing function $G:\rp\to\rp$  is given by the following formula
\[\wt{G}(s):=\sup_{t>0}(s\cdot t-G(t)).\]
If $G$ is a Young function, so is $\wt{G}$. 
If $G$ is an $N$-function, so is $\wt{G}$.

Having Young functions $G,\wt G$, we are equipped with Young's inequality reading as follows\begin{equation}
\label{in:Young} ts\leq G(t)+\wt{G}(s)\quad\text{for all }\ s,t\geq 0.
\end{equation}

 We say that a function $G:\rp\to\rp$ satisfies $\Delta_2$-condition if there exist $c_{\Delta_2},t_0>0$ such that $G(2t)\leq c_{\Delta_2}G(t)$ for $t>t_0.$ 
 We say that $G$ satisfy $\nabla_2$-condition if $\wt{G}\in\Delta_2.$ 
Note that it is possible that $G$ satisfies only one of the conditions $\Delta_2/\nabla_2$. For instance, for $G(t) =( (1+|t|)\log(1+|t|)-|t|)\in\Delta_2$, its complementary function is  $\widetilde{G}(s)= (\exp(|s|)-|s|-1 )\not\in\Delta_2$. See~\cite[Section~2.3, Theorem~3]{rao-ren} for equivalence of various definitions of these conditions and~\cite{CGZG,DFL} for illustrating the subtleties. In particular, $G\in\Delta_2\cap\nabla_2$ if and only if $1<i_G\leq s_G<\infty$, see~\eqref{iG-sG}. This assumption implies a comparison with power-type functions i.e.
$\frac{G(t)}{t^{i_G}}$ 
 is non-decreasing 
  and 
  $\frac{G(t)}{t^{s_G}}$ 
  is non-increasing, but it is stronger than being sandwiched between power functions.
 
 \begin{lem}\label{lem:equivalences}
 If {an $N$-function} $G\in\Delta_2\cap\nabla_2$, then 
 $g(t)t\approx G(t)$ 
  and 
 $\wt G(g(t))\approx G(t)$ 
 with the constants depending only on the growth indexes of $G$, that is $i_G$ and $s_G$.  Moreover,  $g^{-1}(2t)\leq c g^{-1}(t)$ with $c=c(i_G,s_G).$
 \end{lem}

Due to Lemma~\ref{lem:equivalences} and~\cite[Lemmas~3 and~21]{DiEt}, we have the following relations.
\begin{lem}\label{lem:DiEt-mon} If $G$ is an $N$-function of class $C^2((0,\infty))\cap C([0,\infty))$, $G{,g}\in\Delta_2\cap\nabla_2$, $\opA$ is given by~\eqref{opA:def}, then for every $\xi,\eta\in\rnm$ it holds 
\begin{equation}
    \label{opA:strict-monotonicity}
\left(\opA(x,\xi)-\opA(x,\eta)\right):{(\xi-\eta)}\gtrsim_{\data}\, \frac{g(|\xi|+|\eta|)}{|\xi|+|\eta|}|\xi-\eta|^2\approx_{\data} \left| \calV(\xi)-\calV(\xi) \right|^2,
\end{equation}
and\begin{equation}
    \label{relation:g-V} g(|\xi|+|\eta|)|\xi-\eta|\approx_{\data} G^\frac{1}{2}(|\xi|+|\eta|) |\calV(\xi)-\calV(\eta)|.
\end{equation}
\end{lem}
 \subsection{Orlicz spaces}
 Basic reference  for this section is~\cite{adams-fournier}, where the theory of Orlicz spaces is presented for scalar functions. The proofs for functions with values in $\Rm$ can be obtained by obvious modifications.
 
We study the solutions to PDEs in the Orlicz-Sobolev spaces equipped with a~modular function $G\in C^1 {((0,\infty))}$ - a strictly increasing and convex function   such that $G(0)=0$ and satisfying~\eqref{iG-sG}. Let us define a modular \begin{equation}
    \label{modular}
\vr_{G,U}(\texi)=\int_U G(|\texi|)\,dx.
\end{equation}
%
For any bounded $\Omega\subset\rn$, by Orlicz space ${L}^G(\Omega,\Rm)$  we understand the space of measurable functions endowed with the Luxemburg norm 
\[||\tef||_{L^G(\Omega)}=\inf\left\{\lambda>0:\ \ \vr_{G,\Omega}\left( \tfrac{1}{\lambda} |\tef|\right)\leq 1\right\}.\]
 We define the Orlicz-Sobolev space  $W^{1,G}(\Omega)$  as follows
\begin{equation*} 
W^{1,G}(\Omega,\Rm)=\big\{\tef\in {W^{1,1}(\Omega,\Rm)}:\ \ |\tef|,|D{\tef}|\in L^G(\Omega,\Rm)\big\},
\end{equation*}where the gradient is understood in the distributional sense, endowed with the norm
\[
\|\tef\|_{W^{1,G}(\Omega,\Rm)}=\inf\bigg\{\lambda>0 :\ \    \vr_{G,\Omega}\left( \tfrac{1}{\lambda} |\tef|\right)+ \vr_{G,\Omega}\left( \tfrac{1}{\lambda} |D{\tef}|\right)\leq 1\bigg\} 
\]
and  by $W_0^{1,G}(\Omega,\Rm)$ we denote the closure of $C_c^\infty(\Omega,\Rm)$ under the above norm. 
Since condition~\eqref{iG-sG} imposed on $G$ implies $G,\wt{G}\in\Delta_2$, the Orlicz-Sobolev space $W^{1,G}(\Omega,\Rm)$ we deal with is separable and reflexive. Moreover, one can apply arguments of \cite{Gossez} to infer density of  smooth functions  in $W^{1,G}(\Omega,\Rm)$.


The counterpart of the H\"older inequality in this setting reads \begin{equation}
\label{in:Hold} \|\tef\teg\|_{L^1(\Omega,\Rm)}\leq 2\|\tef\|_{L^G(\Omega,\Rm)}\|\teg\|_{L^{\wt{G}}(\Omega,\Rm)}
\end{equation}
for all $\tef\in L^G(\Omega,\Rm)$ and $ \teg\in L^{\wt{G}}(\Omega,\Rm)$.



\subsection {The operator} We notice that in such regime  the operator   $\mathsf{A}_{G}$ acting as
\begin{flalign*}
\langle\mathsf{A}_{G}u,\teph\rangle:=\int_{\Omega}\opA(x,\teDu): \teDph \, dx\quad \text{for}\quad \teph\in C^{\infty}_{0}(\Omega,\Rm)
\end{flalign*}
is well defined on a reflexive and separable Banach space $W^{1,G}(\Omega,\Rm)$ and $\mathsf{A}_{G}(W^{1,G}(\Omega,\Rm))\subset (W^{1,G}(\Omega,\Rm))'$. Indeed, when $\teu\in W^{1,G}(\Omega,\Rm)$ and $\teph\in C_c^\infty(\Omega,\Rm)$, structure condition~\eqref{opA:def},  H\"older's inequality~\eqref{in:Hold}, and Lemma~\ref{lem:equivalences}  justify that
\begin{flalign*}
|{\langle \mathsf{A}_{G}u,\teph \rangle}|\le &\,c\,\int_{\Omega}g(|{\teDu}|){|\teDph|} \, dx \le c\left \| g(|{\teDu}|) \right \|_{L^{\wt G(\cdot)} }\|{|\teDph|}\|_{L^{G}}\nonumber \\
\le &\, c\|{|\teDu|}\|_{L^{ G}}\|{|\teDph|}\|_{L^{G}}\le c\|{\teph}\|_{W^{1,G}}.
\end{flalign*}

\subsection{Definitions of solutions and comments on existence results}
\label{ssec:sols} 

We stress that the problems are considered under the regime of {\rm Assumption {\bf (A-vect)}}.

A function $\tev\in W_{loc}^{1,G}(\Omega,\Rm)$ is called an \underline{$\opA$-harmonic map} in $\Omega\subset\Rn$ provided  \begin{equation}\label{pre-hom-eq}
    \int_\Omega \opA(x,\teDv):\teDvp\,dx=0\qquad\text{for all }\ \tevp\in C_c^\infty(\Omega,\Rm).
\end{equation} 
As a consequence of our main result, in Corollary~\ref{coro:ah-cont}, we prove that $\opA$-harmonic maps are continuous. In fact, by Campanato's characterization \cite[Theorem~2.9]{giusti} one can infer from Proposition~\ref{prop:osc} H\"older continuity  $C^{0,\gamma}(\Omega,\Rm)$ of $\opA$-harmonic maps with any exponent $\gamma\in (0,1).$

A function $\teu\in W^{1,G}_{loc}(\Omega,\Rm)$ is called {a} {weak solution} to~\eqref{eq:mu}, if\begin{equation}
    \label{eq:main-mu-weak}\int_\Omega \opA(x,\teDu): \teDvp\,dx=\int_\Omega\tevp\,d\temu(x)\quad\text{for every }\ \tevp\in C^\infty_c(\Omega,\Rm).
\end{equation}
Recall that $W^{1,G}_{0}(\Omega,\Rm)$ is separable and by its definition ${C_c^\infty}(\Omega,\Rm)$ is dense there.
 
\begin{rem}[Existence and uniqueness of weak solutions] \rm For $\temu\in (W^{1,G}_{0}(\Omega,\Rm))',$ due to the strict monotonicity of the operator, there exists a unique weak solution to~\eqref{eq:mu}, see~\cite[Section~3.1]{KiSt}.
\label{rem:weak-sol}
\end{rem}

Recall that the notion of SOLA is defined in Section~\ref{ssec:main-results}. The problem~\eqref{eq:mu} admits a solution of this type for arbitrary bounded measure.
\begin{prop}\label{prop:exist}
If a vector field $\opA$ satisfies {\rm Assumption {\bf (A-vect)}} and $\temu\in\Mb$, then there exists a SOLA $u\in W^{1,g}(\Omega)$ to~\eqref{eq:mu}.
\end{prop}
The idea to prove it is to 
consider
\[\tef_k(x):= \int_\Rm \vr_k(x-y)\, d\temu(y),\]
where $\vr_k$ stands for a standard mollifier i.e. for a nonnegative, smooth, and even function such that $\int_\Rm\vr(s){\rm\,d}s=1$ we define  $\vr_k(s)=k^n\vr(ks)$ for $k\in \N$. Of course   \[\tef_k\xrightharpoonup[]{*}\temu\qquad\text{and}\qquad\sup_k\|\tef_k\|_{L^1(\Omega)}\leq|\temu|(\rn)<\infty.\]
By Remark~\ref{rem:weak-sol} one finds $\teu_k\in W_0^{1,G}(\Omega,\Rm)$
such that for every $\tevp\in W_0^{1,G}(\Omega,\Rm)$ it holds that
\[\int_\Omega \opA(x,\teDu_k): \teDvp\,dx=\int_\Omega\tevp\,\tef_k dx.\]
The existence of SOLA by passing to the limit can be justified by modification of arguments of \cite{DHM}.

\medskip

\noindent See \cite{ACCZG,CGZG,IC-measure-data,CiMa} for related existence and regularity results in the scalar case and~\cite{DHM,DS,FR,L1,L2} for vectorial existence results in various regimes.

\subsection{Auxiliary results}

\begin{lem}\label{lem:excess} For $\tens{g}:B\to\R^k$, $k\geq 1$ and any $\texi\in\R^k$ it holds that
\[\barint_B|\tens{g}-(\tens{g})_B|\,dx\leq 
2\,\barint_B|\tens{g}-\texi|\,dx.\]
\end{lem}

We have the following corollary of the Cavalieri Principle.
\begin{lem}
\label{lem:cava}
If $\nu\in\mathcal{M}(\Omega)$ has a density $\omega$ (i.e. $d\nu=\omega(x)\,dx$ with $\omega\in L^1(\Omega)$) and $(1+|f|)^{-(\gamma+1)}\omega\in L^1(\rn)$ for some $\gamma>0$, then\[\int_0^\infty \frac{\nu(\{|f|<t\})}{(1+t)^{2+\gamma}}\,dt=\frac{1}{1+\gamma}\int_\rn\frac{d\nu}{(1+|f|)^{\gamma+1}}.\]
\end{lem}

\begin{lem}[\cite{giusti}, Lemma~6.1]
\label{lem:absorb1}
Let $\phi:[R/2,3R/4]\to[0,\infty)$ be a function such that
\[\phi(r_1)\leq \frac{1}{2}\phi(r_2)+\mathsf{A}+\frac{\mathsf{B}}{(r_2-r_1)^\beta}\qquad\text{for every}\qquad R/2\leq r_1<r_2\leq 3{R}/4\]
with $\mathsf{A,B}\geq 0$ and $\beta>0$. Then there exists $c=c(\beta)$, such that\[\phi(R/2)\leq c \left(\mathsf{A}+\frac{\mathsf{B}}{R^\beta}\right).\]
\end{lem}

\begin{lem}[\cite{HanLin}, Lemma~3.4]
\label{lem:absorb2}
Let $\phi(t)$ be a nonnegative and nondecreasing function on $[0,R]$. Suppose that
\[\phi(\rho) \leq \mathsf{A} \left[ \left(\frac{\rho}{r}\right)^{\alpha} + \epsilon \right] \phi(r) + \mathsf{B} r^{\beta}\]
for any $0 < \rho \leq r \leq R$, with $\mathsf{A},\mathsf{B},\alpha,\beta$ nonnegative constants and $\beta < \alpha$.
Then for any $\gamma \in (\beta, \alpha)$, there exists a constant $\epsilon_{0} = \epsilon_{0}(\mathsf{A},\alpha, \beta, \gamma)$ such that if $\epsilon < \epsilon_0$, then for all $0 < \rho \leq r \leq R$ we have
\[ \phi(\rho) \leq c \left \{ \left( \frac{\rho}{r} \right)^{\gamma} \phi(r) + \mathsf{B} \rho^{\beta} \right\}\]
where $c$ is a positive constant depending on $\mathsf{A}, \alpha, \beta, \gamma$.
\end{lem}

The next lemma is a self-improving property for the reverse H\"older inequalities. 

\begin{lem}[\cite{hekima}, Lemma~3.38]\label{lem:self}
Let $0<r<q<p<\infty$, $0< \rho < R \leq 1$ and $w \in L^{p}(B_1)$.
If the following reverse H\"older inequality holds
$$\bigg( \int_{B_{\sigma'}} w^p \, dx \bigg)^{1/p} \leq \frac{c_0}{(\sigma -\sigma')^{\varkappa}} \bigg( \int_{B_{\sigma}} w^{q} \, dx \bigg)^{1/q} + c_1$$
for some constants $c_0$ and $c_1$, whenever $\rho \leq \sigma' \leq \sigma \leq R$.
Then there exists $c=c(c_0,\xi,p,q,r)$ such that
$$\bigg( \int_{B_{\rho}} w^p \, dx \bigg)^{1/p} \leq \frac{c}{(R - \rho)^{\wt{\varkappa}}} \bigg( \int_{B_{R}} w^{r} \, dx \bigg)^{1/r} + c_1,$$
where
$$\wt{\varkappa}=\frac{\varkappa r(p-q)}{q(p-r)}.$$
\end{lem}

The following modular version of the Sobolev-Poincar\'e inequality follows almost directly as in~\cite{Baroni-Riesz}, but we present the proof for vector-valued functions for the sake of completeness.
\begin{prop} \label{prop:OrSob}
Suppose $\Omega$ is a bounded Lipschitz domain in $\rn$, $m,n\geq 1$, and $G:\rp \to\rp$ is an $N$-function such that $G\in\Delta_2\cap\nabla_2$. Then there exist a constant $C=C(n,m,|\Omega|,G)>0$, such that for every $\teu\in W_0^{1,G}(\Omega,\Rm)$  
\[\int_\Omega G^{n'}(|\teu|)\,dx\leq C\left(\int_\Omega G (|\teDu|)\,dx\right)^{n'}.\]
 \end{prop}
\begin{proof} We provide the proof only in the case of continuously differentiable $G$. Otherwise every time one can find a sufficiently smooth function $G_\circ$ comparable to $G$, i.e.~such that there exists $c>0$, such that $G_\circ(t)/c\leq G(t)\leq cG_\circ(t).$ Moreover, we start with the proof for fixed $\teu\in C_c^\infty(\Omega,\Rm)$ and then conclude by the density argument. The classical Sobolev in $W^{1,1}$ inequality gives
 \begin{equation}
 \label{SobN'}
 \left(\int_\Omega G^{n'}(|\teu|)\,dx\right)^\frac{1}{n'}\leq c\,\int_\Omega |{D(G(| \teu|))}|\,dx .
 \end{equation}
Since $G\in\Delta_2,$ it satisfies $g(t)\leq c\, G(t)/t$ and $G^*\left(\frac{G(t)}{t}\right)\leq {G(t)}.$ Thus by the Young inequality  we arrive at
 \begin{equation}
 \label{grad-vp-est}\begin{split}
 |{D(G(| \teu|))}|&= g(|\teu|)\big|{D|\teu|}\big|\leq  c  \frac{G(|\teu|)}{\teu}|\teDu|\\
 &\leq\ve G^*\left(\frac{G(|\teu|)}{|\teu|}\right)+c \, G(|\teDu|)\leq \ve G(|\teu|)+c\, G(|\teDu|).\end{split}
 \end{equation}
 Summing up, we have
 \[\left(\int_\Omega G^{n'}(|\teu|)\,dx\right)^\frac{1}{n'}\leq  C\ve \int_\Omega G(|  \teu|) \,dx+Cc_\ve\int_\Omega G(|\teDu|)\,dx,\]
where according to the H\"older inequality we obtain
 \[\left(\int_\Omega G^{n'}(|\teu|)\,dx\right)^\frac{1}{n'} \leq \ve C|\Omega|^\frac{1}{n}\left( \int_\Omega G^{n'} (|  \teu|) \,dx\right)^\frac{1}{n'}+Cc_\ve\int_\Omega G(|\teDu|)\,dx.\]
 Now we can choose $\ve$ small enough to absorb it on the right-hand side and obtain the claim for $\teu\in C_c^\infty(\Omega,\Rm)$. Since smooth function are dense in $W_0^{1,G}(\Omega,\Rm)$ by standard approximation argument we get the claim for all $\teu\in W_0^{1,G}(\Omega,\Rm)$.
\end{proof}

\subsection{Properties of $\opA$-harmonic maps}

Let us establish some fundamental properties of $\opA$-harmonic maps. 
\begin{prop}[Caccioppoli estimate]\label{prop:Cacc}
If $\tev \in W_0^{1,G} (\Omega,\Rm)$ is a nonnegative $\opA$-harmonic map, $\telambda\in\Rm$, and $7/8\leq\sigma'<\sigma\leq 1,$ then there exists $c=c(\data)>0$, such that  
\begin{equation}
\label{inq:caccI} 
\barint_{B_{\sigma' r}}G(|\teDv|)\,dx\leq\frac{c}{(\sigma'-\sigma)^{s_G}}\barint_{B_{\sigma r}}G\left(\frac{|\tev-\telambda|}{r}\right)
\,dx.
\end{equation}
\end{prop}
\begin{proof}  Let us pick a~cutoff function $\eta\in C_c^\infty({B_{\sigma r}},\Rm)$ such that $\mathds{1}_{B_{\sigma' r}}\leq\eta\leq \mathds{1}_{B_{\sigma r}}$ and $|\teDet|\leq c_1/(\sigma'-\sigma)$. We use $\xi=\eta^q (v-\telambda)$ as a test function to get
\[\int_\Omega \opA(x,\teDv): \teDv\,\eta^q \,dx =\int_\Omega \opA(x,\teDv):(-q\eta^{q-1}(\tev-\telambda)\otimes \teDet)\,  \,dx
.\]
Therefore, due to coercivity of $\opA$ and the Cauchy-Schwartz inequality we have
\begin{equation*}\int_{B_{\sigma r}}G(|\teDv|)\eta^q\,dx\leq c\,\int_{B_{\sigma r}}\frac{g(|\teDv|)}{|\teDv|}\eta^{q-1}|\teDv:((\tev-\telambda)\otimes\teDet)| \,dx=:\mathcal{K}
\end{equation*}
Noting that $q$ is large enough to satisfy $s_G'\geq q'$, we have in turn that $\wt{ G}(\eta^{q-1}t)\leq c \eta^q\wt{G}(t)$  and  \[\wt{G}(\eta^{q-1}g(t))\leq c \eta^q \wt{G}(g(t))\leq c \eta^q G(t).\] 
Then,  using Young inequality~\eqref{in:Young} applied to the integrand of $\mathcal{K}$ we get
\begin{equation*}
\begin{split}\mathcal{K}&\leq \ve \int_{B_{\sigma r}}\wt{G}(\eta^{q-1}|\teDv|)\, dx+ c_\ve \int_{B_{\sigma r}}G\left(|\tev-\telambda|\,|\teDet| /c_1 \right)\,
dx\\& \leq \ve c\,\int_{B_{\sigma r}} \eta^{q }G(|\teDv| ) \,dx+ c_\ve \int_{B_{\sigma r}} G\left(|\tev-\telambda|\,|\teDet| /c_1 \right)\,
dx
\end{split}
\end{equation*}
with arbitrary $\ve<1$. Choosing $\ve$ small enough to absorb the term, and finally by properties of $\eta$ we obtain~\eqref{inq:caccI}.
\end{proof}

An $\opA$-harmonic function $\tev$ is a minimizer of a functional 
\[ \tev \mapsto \int_{\Omega} G(|\teDv|) \, dx. \]
Therefore, taking the operator independent of the spacial variable\begin{equation}
    \label{opA-no-x}
\opA(x,\xi)=c_a\frac{g(|\xi|)}{|\xi|}\xi\end{equation}
$\opA$-harmonic functions are Lipschitz regular by the following fact.
\begin{lem}[\cite{DSV1}, Lemma~5.8] Suppose $G$ is an $N$-function of class $C^2((0,\infty))\cap C([0,\infty))$, $G,g\in\Delta_2\cap\nabla_2$, $\tew$ is $\opA$-harmonic in $\Omega$ for $\opA$ given by~\eqref{opA-no-x}. Let $B\subset 2B\Subset\Omega$.  \label{lem:DSV1} Then there exists $c=c(\data)>0$ such that \[\sup_{B} G(|D\tew|) \leq c\, \barint_{2B} G(|D\tew|) \, dx.\]
\end{lem}

\begin{prop} If $\tev\in W^{1,G}(\Omega,\Rm)$ is $\opA$-harmonic, then for any $\varsigma \in (0,1)$ there exists $R_0=R_0(\data,\varsigma) \in (0,1]$ such that for any $0< R \leq R_0$ and $B_R\subset \subset \Omega$ it holds that 
\label{prop:osc}
\begin{equation}\label{osc-dec-est}
\barint_{B_{\delta R}} |\tev - (\tev)_{B_{\delta R}}| \, dx \leq c_{\rm o} \delta^{1+(\varsigma - 1)/{s_{G}}} \barint_{B_{R}} |\tev - (\tev)_{B_{R}}| \, dx
\end{equation}
whenever $\delta \in (0,1/4],$ where $c_{\rm o}\geq 1$ is a constant depending only on $\data$ and $\varsigma$.
\end{prop}
\begin{proof} 
In this proof, we use a classical perturbation argument, see for instance \cite[Theorem 3.8]{HanLin}.
It suffices to prove \eqref{osc-dec-est} for $\tev$ solving \[-\tedv\left(a(x)\frac{g(|\teDv|)}{|\teDv|}\teDv\right)=0\quad\text{ in }\quad B_1.\] Indeed, the general case can be deduced then by considering $\wt{\tev}(x)=
{\tev(x_0+Rx)}/{R}$ solving 
$-\tedv\bopA(x,{D\wt{\tev}})=0$ on $B_1(0)$ with \[\bopA(x,\xi)=\bar{a}(x)\frac{\bar g(|\xi|)}{|\xi|}\xi=\frac{1}{g(R)}\opA(x_0+Rx,R\xi)= {a(x_0+Rx)}\frac{g(R|\xi|)}{g(R)|\xi|}{\xi}.\] In this case, the modulus $\omega_{\bar a}$ of continuity of $\bar a$ satisfies $\omega_{\bar a}(r) = \omega_{a}(rR)$.

Note first Propositions~\ref{prop:OrSob} and~\ref{prop:Cacc} imply that for any $7/8 \leq \sigma' < \sigma \leq 1$
\[ \left( \barint_{B_{\sigma'}} G^{n'}\left( \frac{|\tev - (\tev)_{B_{\sigma'}}|}{\sigma'} \right) \, dx \right)^{1/n'} \leq \frac{c}{(\sigma-\sigma')^{s_G}} \barint_{B_{\sigma}} G\left( \frac{|\tev-(\tev)_{B_{\sigma}}|}{\sigma} \right) \, dx.\]
From the doubling property of $G$ and the upper and lower bound on $\sigma'$ and $\sigma$, we have
\[ \left( \barint_{B_{\sigma'}} G^{n'}\left( |\tev - (\tev)_{B_{\sigma'}}| \right) \, dx \right)^{1/n'} \leq \frac{c}{(\sigma-\sigma')^{s_G}} \barint_{B_{\sigma}} G\left( |\tev-(\tev)_{B_{\sigma}}| \right) \, dx. \]
Using the triangular inequality and Jensen's inequality
\begin{align*}
& \left( \barint_{B_{\sigma'}} G^{n'}\left( |\tev - (\tev)_{B_{1}}| \right) \, dx \right)^{1/n'} \\
&\qquad\qquad \leq c \left( \barint_{B_{\sigma'}} G^{n'}\left( |\tev - (\tev)_{B_{\sigma'}}| \right) \, dx \right)^{1/n'} + c\,G \left( |(\tev)_{B_{\sigma'}} - (\tev)_{B_{1}}| \right) \\
&\qquad\qquad \leq \frac{c}{(\sigma-\sigma')^{s_G}} \barint_{B_{\sigma}} G\left( |\tev-(\tev)_{B_{\sigma}}| \right) \, dx + c\,\barint_{B_{\sigma'}}G\left( |\tev - (\tev)_{B_{1}}| \right) \, dx \\
&\qquad\qquad \leq c\left(\frac{1}{(\sigma-\sigma')^{s_G}} + 1 \right) \barint_{B_{\sigma}} G\left( |\tev-(\tev)_{B_{1}}| \right) \, dx + c\,G \left( |(\tev)_{B_{1}} - (\tev)_{B_{\sigma}}| \right) \\
&\qquad\qquad \leq  \frac{c}{(\sigma-\sigma')^{s_G}} \barint_{B_{\sigma}} G\left( |\tev-(\tev)_{B_{1}}| \right) \, dx,
\end{align*} for $c=c(\data)$. From Lemma~\ref{lem:self}, for any $t \in (0,1/s_{G})$, we have
\begin{align} \label{self-imp-vexc}
\left( \barint_{B_{7/8}} G^{n'}\left(|\tev - (\tev)_{B_{7/8}}|\right) \, dx \right)^{1/n'}
& \leq c \left( \barint_{B_{7/8}} G^{n'}\left(|\tev - (\tev)_{B_{1}}|\right) \, dx \right)^{1/n'} \notag\\
& \leq c\left( \barint_{B_{1}} G^{t}\left(|\tev - (\tev)_{B_{1}}|\right) \, dx \right)^{1/t} \notag \\
& \leq c \,G \left( \barint_{B_{1}} |\tev - (\tev)_{B_{1}}| \, dx \right).
\end{align}
Here, in the last line, as $t \mapsto G^{t}(t)$ is a concave function for every $t \in (0, 1/s_{G})$ we have used Jensen's inequality.

Let $\tew \in \tev + W^{1,G}_{0}(B_{\sigma}, \R^m)$ for any $\sigma\in(1,1/2)$ be the weak solution to
\begin{equation}\label{pre-ref-eq}
- \tedv \left(a(x_0) \frac{g(|\teDw|)}{|\teDw|}\teDw\right) = 0 \qquad \text{in } B_{\sigma}.
\end{equation}
Recall that the function $\calV$ describing ellipticity of $\opA$ is defined in~\eqref{V-def} and $a$ is a function bounded below by ${c_a}>0$ and with a modulus of continuity $\omega_a$. Testing~\eqref{pre-ref-eq} and~\eqref{pre-hom-eq} against $(\tew - \tev)$, and applying Lemma~\ref{lem:DiEt-mon}, for any $\epsilon \in (0,1]$ we see
\begin{align*}
& \int_{B_{\sigma}} c_a \left| \calV(\teDw)-\calV(\teDv) \right|^2 \, dx \\
&\qquad\qquad\leq c \int_{B_{\sigma}} a(x_0) \left(\frac{g(|\teDw|)}{|\teDw|} \teDw- \frac{g(|\teDv|)}{|\teDv|} \teDv\right) :  (\teDw - \teDv) \, dx \\
&\qquad\qquad = c \int_{B_{\sigma}} (a(x_0) - a(x)) \frac{g(|\teDv|)}{|\teDv|} \teDv : (\teDw-\teDv) \, dx \\
&\qquad\qquad \leq c \,\omega_{a}(\sigma)\int_{B_{\sigma}} g(|\teDw| + |\teDv|)\, |\teDw-\teDv| \, dx \\
&\qquad\qquad \leq c\,  \epsilon \int_{B_{\sigma}} |\calV(\teDw)-\calV(\teDv)|^2 \, dx + c\,\frac{\omega_{a}(1/2)^2}{\epsilon} \int_{B_{\sigma}} G(|\teDw|+|\teDv|) \, dx,
\end{align*}
where in the last line we used Young's inequality. As $c=c(\data)$ and $c_a>0$, we can take $\epsilon$ small enough to absorb the first term on the right-hand side. Then Jensen's inequality implies
\begin{flalign}\nonumber
\int_{B_{\sigma}} \left| \calV(\teDw)-\calV(\teDv) \right|^2 \, dx &\leq c\, \omega_{a}(1/2)^2 \int_{B_{\sigma}} [G(|\teDv|) + G(|\teDw|)] \, dx.\end{flalign}
Since $\tew$ is a minimizer of the integral functional
\[ \tew \mapsto \int_{B_{\sigma}} G(|\teDw|) \, dx, \]
we have
\begin{flalign}\label{pre-excess-1} 
\int_{B_{\sigma}} \left| \calV(\teDw)-\calV(\teDv) \right|^2 \, dx \leq c\, \omega_{a}(1/2)^2 \int_{B_{\sigma}} G(|\teDv|) \, dx.
\end{flalign}
Since we know the Lipschitz regularity of $\tew$ provided by Lemma~\ref{lem:DSV1} \[\sup_{B_{\sigma/2}} G(|D\tew|) \leq c \,\barint_{B_{\sigma}} G(|D\tew|) \, dx,\]
it follows from \eqref{pre-excess-1} that 
\begin{align*}
\int_{B_{\delta \sigma}} G(|\teDv|) \, dx &\leq c\,
\int_{B_{\delta \sigma}} |\calV(\teDv)|^2 \, dx \\
& \leq c\, \int_{B_{\delta \sigma}} |\calV(\teDv) - \calV(\teDw)|^2 \, dx + \int_{B_{\delta \sigma}} |\calV(\teDw)|^2 \, dx \\
& \leq c\, \int_{B_{\sigma}} |\calV(\teDv) - \calV(\teDw)|^2 \, dx + \delta^{n} \int_{B_{\sigma}} |\calV(\teDw)|^2 \, dx \\
& \leq c  \left( \delta^{n} + \omega_{a}(1/2)^2 \right) \int_{B_{\sigma}} G(|\teDv|) \, dx.
\end{align*}
We take $R_0=R_{0}(\data,\varsigma)$ small enough to ensure that $\omega_{a}(R_0)^2<\epsilon_{0}$, where $\epsilon_0$ is a constant given in Lemma~\ref{lem:absorb2}.
Then Lemma~\ref{lem:absorb2}, Proposition~\ref{prop:Cacc} and \eqref{self-imp-vexc} give
\begin{equation}\label{pre-excess-2}
\barint_{B_{\delta}} G(|\teDv|) \, dx \leq c \delta^{\varsigma - 1} \barint_{B_{1/2}} G(|\teDv|) \, dx \leq c \delta^{\varsigma - 1} G \left( \int_{B_1} |\tev - (\tev)_{B_1}| \, dx \right)
\end{equation}
where $c$ depends only on $\data$ and $\varsigma$.

Using the Sobolev-Poincar\'e inequality in $W^{1,1}$ and Jensen's inequality, we conclude that 
\begin{align*}
\barint_{B_{\delta}} |\tev - (\tev)_{B_{\delta}}| \, dx
& \leq c \,\delta \, \barint_{B_{\delta}} |\teDv| \, dx \\
& \leq c \,\delta \, G^{-1} \left( \barint_{B_{\delta}} G(|\teDv|) \, dx \right) \\
& \leq c \, \delta^{1+(\varsigma - 1)/i_{G}} \int_{B_1} |\tev - (\tev)_{B_1}| \, dx,
\end{align*}
what completes the proof.
\end{proof}

\section{Measure data $\opA$-harmonic approximation}\label{sec:Ah-approx}

In this section we provide the tool of crucial meaning for our further reasoning -- the  approximation  of a $W^{1,G}$-function by an $\opA$-harmonic map for weighted operator $\opA$ of an Orlicz growth given by \eqref{opA:def}. Results in this spirit can be found in \cite{DuGr,DuMi2004,DSV3}, but most preeminently for the approximation relevant for application to measure data problems we refer to \cite[Theorem~4.1]{KuMi2018}.

We define an auxiliary function\begin{equation}\label{Hs-def}
    H_s(t)=\int_{0}^{t} \frac{g(r)^{1-s} G(r)^{s}}{r} \, dr\qquad\text{for }\ s \in [0,1/2).
\end{equation}
It is readily checked that when $s>\max\{2-i_{G},0\}$, $H_s$ is a Young function satisfying \begin{equation}
    \label{Hs-prop-1}
H_s(t) \approx g(t)^{1-s} G(t)^{s}
\end{equation} with intrinsic constants depending only on $i_G,s_G$ and $s$. In fact, $H_s\in\Delta_2\cap\nabla_2$ since \begin{equation}
    \label{Hs-prop-2} 0 < s + i_G -2\leq \frac{t H_s''(t)}{H_s'(t)} = \frac{(1-s) t g'(t)}{g(t)} + \frac{s t g(t)}{G(t)} \leq s + s_G -2.\end{equation}
Furthermore, there exist $\epsilon,c,t_0>0$, such that $H_s(t)\geq c t^{1+\epsilon}$ and $H_s(t)\geq c g^{1+\epsilon}(t)$ for all $t\geq t_0.$

\begin{theo}\label{theo:Ah-approx} Under Assumption {\bf (A-vect)} let $\ve>0$, $\gamma\in (0,1/(s_Gn))$, and \begin{equation}
    \label{s-range} 
    \max\{2-i_{G},0\} < s  <s_{\rm m}:=\frac{i_G-\gamma s_G n}{i_G+ s_Gn}.
\end{equation}Suppose that $\teu\in W^{1,G}(B_r(x_0),\Rm)$ satisfies\begin{equation}
    \label{u-male}\barint_{B_r(x_0)}|\teu|\,dx\leq Mr,\quad M\geq 1
\end{equation}
then there exists $\delta=\delta(\data,s,M,\ve)\in(0,1]$ such that if $\teu$ is almost $\opA$-harmonic in a sense that 
for every $\tevp\in W^{1,G}_0(B_r(x_0),\Rm)\cap L^\infty(B_r(x_0),\Rm)$ it holds 
\begin{equation}
    \label{Du-male}\left|\barint_{B_r(x_0)}  
    \opA(x,\teDu):\teDvp\,dx\right|\leq\frac{\delta}{r}\|\tevp\|_{L^\infty(B_r(x_0),\Rm)},
\end{equation} 
then there exists an $\opA$-harmonic map $\tev\in W^{1,G}(B_{r/2}(x_0),\Rm)$ satisfying
    \begin{equation}\label{Du-blisko-Dv}
    \barint_{B_{r/2}(x_0)}H_s(|\teDu-\teDv|)\,dx
    \leq \ve
\end{equation}
together with\begin{equation}\label{v-male}
    \barint_{B_{r/2}(x_0)}|\tev|\,dx\leq 2^nMr
    \quad\text{and}\quad 
    \barint_{B_{r/2}(x_0)}H_s(|\teDv|)\,dx \leq cH_s(M),
\end{equation} 
where $c=c(\data)>0.$
\end{theo}
\begin{rem} \label{rem:iG-male} The limitation that $G$ has to be superquadratic ($i_G\geq 2$) can be a~little bit relaxed in Theorem~\ref{theo:Ah-approx} and later on the restriction is not needed.  The key property is to ensure that the range of admissible $s$ from  \eqref{s-range} is {nonempty.  We}  need to assume that $i_G$ is either bigger or equal to $2$, or close to $2$ in a sense that
\[2-i_G <\frac{i_G-1}{i_G+ s_Gn}.\]
\end{rem}
\begin{proof}The plan is to first establish suitable a priori estimates for the rescaled problem and then proceed with the proof via contradiction. The proof is presented in $6$ steps.

\textbf{Step 1. Scaling. } We fix arbitrary $\tevp\in W^{1,G}_0(B_r(x_0),\Rm)\cap L^\infty(B_r(x_0),\Rm)$ satisfying~\eqref{Du-male}.  Let us change variables putting\begin{equation}
    \label{bar}
    \tebu(x):=\frac{\teu(x_0+rx)}{Mr},\  \bopA(x_0+rx,\texi)=\opA(x_0+rx,M\texi),\ \text{and}\  \teeta(x):=\frac{\tevp(x_0+rx)}{r}.
\end{equation}
Then $\bopA$ satisfies the same conditions as $\opA$ with the functions $
\bg(t):=g(Mt)$ and 
$\bG(t):=G(Mt)/M$ 
(with $\bG'=\bg$), and the constants depending on $\data.$ Of course in such a case $i_G=i_\bG$ and $s_G=s_\bG.$

Having~\eqref{u-male} and~\eqref{Du-male}, by denoting the unit ball by $B_1$, we get further that
\begin{flalign}
\label{bu-male}&\barint_{B_1}|\tebu|\,dx\leq  1,\\
\label{Dbu-male} &\left|\barint_{B_1} \bopA(x,\teDbu):\teDet\,dx\right|\leq {\delta}\|\teeta\|_{L^\infty(B_1,\Rm)}.
\end{flalign}

\textbf{Step 2. A priori estimates. } We choose $q\geq s_G$, pick 
\begin{equation}
    \label{eta}\teeta:=\phi^qT_k(\tebu)\quad\text{with some}\ \phi\in C_c^\infty(B_1),\ 0\leq\phi\leq 1,\ k\geq 0,
\end{equation}
and denote\[\P:=\frac{\tebu\otimes\tebu}{|\tebu|^2}.\]
Then
\begin{flalign*}
\teDet=&\mathds{1}_{\{|\tebu|\leq k\}}(\phi^q \teDbu+q\phi^{q-1}\tebu\otimes D\phi)\\
&+\mathds{1}_{\{|\tebu|> k\}}(\phi^q ({\id}-{\P}) \teDbu+q\phi^{q-1}\tebu\otimes D\phi).
\end{flalign*}
We use~\eqref{eta} in~\eqref{Dbu-male} to get\begin{flalign}
\nonumber& \left|\int_{B_1\cap \{|\tebu|\leq k\}} \bopA(x,\teDbu):(\phi^q \teDbu+q\phi^{q-1}\tebu\otimes D\phi)\,dx\right.\\
&\left.+\int_{B_1\cap \{|\tebu|> k\}} \bopA(x,\teDbu):(\phi^q (\id-\P) \teDbu+q\phi^{q-1}\tebu\otimes D\phi)\,dx\right|\leq {\delta}|B_1|\|\teeta\|_{L^\infty(B_1,\Rm)}.\label{est-1}
\end{flalign}
Since $\bopA$ has the quasi-diagonal structure resulting from~\eqref{opA:def} and
\[  \teDbu:\big((\id-\P)\teDbu\big)=|\teDbu|^2-\frac{D_j\tebu^\alpha\tebu^\alpha D_j\tebu^\beta\tebu^\beta}{|\tebu|^2}=|\teDbu|^2-\frac{\sum_{j=1}^m\langle D_j \tebu ,\tebu\rangle ^2}{|\tebu|^2}\geq 0,
\]
we infer that 
\begin{equation}\label{calc-eta}\bopA(x, \teDbu):\big((\id-\P)\teDbu\big)\geq 0.
\end{equation}

By rearranging terms in~\eqref{est-1}, applying {Lemma~\ref{lem:equivalences},} noting that $\|\teeta\|_{L^\infty(B_1,\Rm)}\leq k$, and dropping a~nonnegative term due to~\eqref{calc-eta}, we get for some $c=c(\data,q)$\begin{flalign}
\nonumber \int_{B_1\cap \{|\tebu|\leq k\}}\bG(|\teDbu|) &\phi^q\,dx\leq c \int_{B_1\cap \{|\tebu|\leq k\}} \frac{\bg(|\teDbu|)}{|\teDbu|}\phi^{q-1}\big|\teDbu:(\tebu\otimes D\phi)\big|\,dx\\
+\,c&\int_{B_1\cap \{|\tebu|> k\}} \frac{k}{|\tebu|} \frac{\bg(|\teDbu|)}{|\teDbu|}\phi^{q-1}\big|\teDbu:(\tebu\otimes D\phi)\big|\,dx+ c|B_1|\delta k.\nonumber 
\end{flalign} 
We estimate the first term on the right-hand side of the last display by the use of Young inequality with a parameter, use Lemma~\ref{lem:equivalences}, and we absorb one term. The second term can be estimated by Schwartz inequality. Altogether we obtain  \begin{flalign}\nonumber
  \int_{B_1\cap \{|\tebu|< k\}}\bG(|\teDbu|) \phi^q\,dx&\leq c \int_{B_1\cap \{|\tebu|< k\}}\bG(|\tebu|\,|D\phi|) \,dx+ c|B_1|\delta k\\
  &+\,ck \int_{B_1\cap \{|\tebu|\geq k\}} {\bg(|\teDbu|)}\phi^{q-1}| D\phi|\,dx \label{st2-apriori}
\end{flalign}
for some $c=c(\data,q)$. 

\medskip

\textbf{Step 3. Summability of $\tebu$ and $D\tebu$. }  We choose $k=t$ in~\eqref{st2-apriori}, then multiply this line by $(1+t)^{-(\gamma+2)},$ where $\gamma>0,$ integrate it from zero to infinity and apply Cavalieri's principle (Lemma~\ref{lem:cava}) twice (with $\nu_1=\bG(|\teDbu|) \phi^q$ and $\nu_2=\bG(|\tebu|\,|D\phi|)$. Altogether we get
 \begin{flalign}\nonumber
 \frac{1}{1+\gamma} \int_{B_1 }\frac{\bG(|\teDbu|) \phi^q}{(1+|\tebu|)^{1+\gamma}}\,dx&\leq
 \frac{c}{1+\gamma} \int_{B_1 }\frac{\bG(|\tebu|\,|D\phi|) }{(1+|\tebu|)^{1+\gamma}}\,dx+ \frac{c}{\gamma}\delta\\
  &+\, c\,\int_0^\infty  \frac{t}{(1+t)^{\gamma+2}}\int_{B_1\cap \{|\tebu|> t\}} {\bg(|\teDbu|)}\phi^{q-1}| D\phi|\,dx\,dt.\label{post-cava}
\end{flalign}
The right-most term in the last display can be estimated as follows
\begin{flalign}\nonumber
c\,\int_0^\infty \frac{t}{(1+t)^{\gamma+2}}\int_{B_1\cap \{|\tebu|> t\}}& {\bg(|\teDbu|)}\phi^{q-1}| D\phi|\,dx\,dt\\\nonumber
&\leq c \int_0^\infty \frac{1}{(1+t)^{\gamma+1}}dt\,\int_{B_1 } {\bg(|\teDbu|)}\phi^{q-1}|D\phi|\,dx \\
&\leq \frac{c}{\gamma}  \int_{B_1 } {\bg(|\teDbu|)}\phi^{q-1}| D\phi|\,dx .\label{do-younga}
\end{flalign}
To estimate it further we note that $q$ is large enough to satisfy $s_G'\geq q'$, {there exist $c_0,c_1>0$ depending on $i_G,s_G$, such that we have ${\bG}^*(c_0\phi^{q-1}\bg(t))\leq c_1 \phi^q {\bG}^*(\bg(t))\leq \phi^q \bG(t).$} Then,  using Young inequality~\eqref{in:Young} applied to the integrand in \eqref{do-younga} and the above observation we get
\begin{equation*}
\begin{split}\int_{B_1 } &{\bg(|\teDbu|)}\phi^{q-1}|D\phi|\,dx\\
&\leq  \frac{1}{2(1+\gamma)}  \int_{B_{1}}\frac{{\bG^*}\left( {c_0}{\phi^{q-1}\bg(|\teDbu|) }\right)}{(1+|\tebu|)^{1+\gamma}} \, dx + \wt c \int_{B_{1}}\frac{\bG\left((1+|\tebu|)^{1+\gamma}|D\phi|\right)}{(1+|\tebu|)^{1+\gamma}}\,dx \\
& \leq \frac{1}{2(1+\gamma)} \int_{B_{1}}\frac{\bG\left(|\teDbu| \right)\phi^{q }}{(1+|\tebu|)^{1+\gamma}}  \,dx + \wt c \int_{B_{1}}\frac{\bG\left((1+|\tebu|)^{1+\gamma}|D\phi|\right)}{(1+|\tebu|)^{1+\gamma}}\, dx
\end{split}
\end{equation*}
with $\wt c=\wt c(\gamma,i_G,s_G)$. By applying this estimate in~\eqref{post-cava} and rearranging terms we obtain  
 \begin{flalign}
 \int_{B_1 }\frac{\bG(|\teDbu|) \phi^q}{(1+|\tebu|)^{1+\gamma}}\,dx&\leq
 {c} \int_{B_{1}}\frac{\bG\left((1+|\tebu|)^{1+\gamma}|D\phi|\right)}{(1+|\tebu|)^{1+\gamma}}\, dx + c\,\delta\frac{1+\gamma}{\gamma} .\label{post-cava-2}
\end{flalign}
 
 {Observe that $1/(s_G n)< i_G-1$, as otherwise the condition required by Remark~\ref{rem:iG-male} is violated}. Recall that since $\gamma<1/(s_G n)$, we have $
\gamma< i_G-1$. Let us set\begin{equation}
    \vt(x):=\frac{\bG((1+|\tebu|)\phi^q)}{(1+|\tebu|)^{1+\gamma}}
\label{vartheta}
\end{equation}
and notice that since $|{D|\tebu|}|\leq|\teDbu|$, by Jensen's inequality, and Lemma~\ref{lem:equivalences} we can estimate 
\begin{flalign}
\nonumber \left|{D\vt}\right|&=\frac{\left|{D( \bG((1+|\tebu|)\phi^q)}(1+|\tebu|)^{1+\gamma}-{D\left((1+|\tebu|)^{1+\gamma}\right)}G((1+|\tebu|)\phi^q)\right|}{(1+|\tebu|)^{2+2\gamma}}\\\nonumber
&\leq c \, \frac{\bg((1+|\tebu|)\phi^q)\left|{D((1+|\tebu|)\phi^q)}\right|}{(1+|\tebu|)^{1+\gamma}} + {(1+\gamma)} \frac{\left|D{\tebu} \right|\bG((1+|\tebu|)\phi^q)} {(1+|\tebu|)^{2+\gamma}}\\\nonumber
&\leq c \, \frac{\bg((1+|\tebu|)\phi^q)|{D\tebu}|\phi^q}{(1+|\tebu|)^{1+\gamma}}\\ \nonumber
&\quad + c \, \frac{\bg((1+|\tebu|)\phi^q)(1+|\tebu|)\phi^{q-1}\left|{D\phi}\right|}{(1+|\tebu|)^{1+\gamma}} + {(1+\gamma)} \frac{  \left|{D \tebu}\right|} {1+|\tebu|}\vt\\
&{\leq  c(\gamma)\frac{  \left|{D \tebu}\right|} {1+|\tebu|}\vt+c  \, \frac{\left|{D\phi}\right|}{\phi}\vt
.}
\label{Dvt}
\end{flalign}
 For the later use in \textbf{Step 6},  we emphasize the dependence of constants on~$\gamma$ by denoting $c(\gamma)$. Note that every $c(\gamma)$ in \eqref{Dvt}-\eqref{int-vt} is an increasing function of~$\gamma$.

Since $q \geq s_G$, for any $\kappa \in [1,i_{G})$ we see
\begin{align*}
\left|{D\left(\vt^{\frac{1}{\kappa}}\right)}\right|
\leq \frac{1}{\kappa} \vt^{\frac{1}{\kappa}-1} \left|{D\vt}\right| 
\leq {c(\gamma)} \, \vt^{\frac{1}{\kappa}} \frac{\left|{D\tebu}\right|}{1+|\tebu|} + {c }\, \vt^{\frac{1}{\kappa}} \frac{\left|{D \phi} \right|}{\phi}
\end{align*}
and 
\begin{align}\label{high-vt}
\left|{D\left(\vt^{\frac{1}{\kappa}}\right)}\right|^{\kappa} 
&\leq c(\gamma) \frac{G((1+|\tebu|)\phi^{q})}{\left[(1+|\tebu|) \phi^{q} \right]^{\kappa}} \frac{\left(|\teDbu | \phi^{q} \right)^{\kappa}}{(1+|\tebu|)^{1+\gamma}}
+ c \frac{G(1+|\tebu|)}{(1+|\tebu|)^{1+\gamma}} |{D \phi}|^{\kappa}.
\end{align}

To proceed further, we define an auxiliary function $h_\kappa$ by setting
$$h_{\kappa}^{-1}(t):= \int_{0}^{t} \frac{1}{\left[G^{-1}(\tau)\right]^{\kappa}} \, d\tau.$$
A straightforward calculation gives
\begin{equation}\label{sob-aux-1}
-1 < - \frac{\kappa}{i_G} \leq \frac{t [h_{\kappa}^{-1}]''(t)}{[h_{\kappa}^{-1}]'(t)}
= - \frac{\kappa t [G^{-1}]'(t)}{G^{-1}(t)} \leq - \frac{\kappa}{s_G} <0,
\end{equation}
which implies that $h_{\kappa}^{-1}$ is an increasing concave function on $[0, \infty)$.
Note that \eqref{sob-aux-1} also gives
$$h_{\kappa}\left( \frac{G(t)}{t^{\kappa}} \right) \approx G(t)$$
with intrinsic constants depending on $i_G,s_G$ and $\kappa$ only. Moreover, we have
$$\bg(t)[h_{\kappa}^{-1}]'(\bG(t))=[h_{\kappa}^{-1}(\bG(t))]' = \frac{d}{dt} \left( \int_{0}^{G(t)} \frac{1}{[G^{-1}(\tau)]^{\kappa}} \, d \tau \right) = \frac{\bg(t)}{t^{\kappa}},$$
and so
\begin{align*}
[h_{\kappa}^{*}]^{-1}(G(t)) \approx h_{\kappa}'(h_{\kappa}^{-1}(\bG(t)) = \frac{1}{[h_{\kappa}^{-1}]'(h_{\kappa}(h_{\kappa}^{-1}(G(t))))}= \frac{1}{[h_{\kappa}^{-1}]'(G(t))} = t^{\kappa}.
\end{align*}
Hence, $h_{\kappa}^{*}(t) \approx \bG(t^{1/\kappa})$.

Applying Young's inequality~\eqref{in:Young} with the pair of  Young functions $(h_{\kappa},h_{\kappa}^{*})$ to \eqref{high-vt}, for any $\epsilon_0 \in (0,1)$ we discover
\begin{align*}
\left|{D\left(\vt^{\frac{1}{\kappa}}\right)}\right|^{\kappa}
& \leq \epsilon_{0} \frac{G((1+|\tebu|)\phi^{q})}{(1+|\tebu|)^{1+\gamma}} + c(\epsilon_0)c(\gamma) \frac{\bG\left(|\teDbu | \phi^{q} \right)}{(1+|\tebu|)^{1+\gamma}}
+ c \frac{G(1+|\tebu|)}{(1+|\tebu|)^{1+\gamma}} |{D \phi}|^{\kappa}.
\end{align*}
By the classical Sobolev inequality we get that\begin{equation} 
    \label{Sob}
   \left( \int_{B_1} |\vt|^\frac{n}{n-\kappa} \,dx\right)^\frac{n-\kappa}{n}\leq c
    \int_{B_1}  \left|{D \left(\vt^{\frac{1}{\kappa}}\right)}\right|^{\kappa} \,dx.
\end{equation}
Merging \eqref{post-cava}, \eqref{post-cava-2}, \eqref{vartheta}, \eqref{Dvt} and \eqref{Sob} and taking $\epsilon_0$ small enough, we get
\begin{equation}\label{int-vt}
    \left(\int_{B_1}\left( \frac{\bG((1+|\tebu|)\phi^q)}{(1+|\tebu|)^{1+\gamma}}\right)^\frac{n}{n-\kappa}\,dx\right)^\frac{n-\kappa}{n}\leq
 \bar c \left( \int_{B_{1}}\frac{\bG\left((1+|\tebu|)^{1+\gamma}|D\phi|\right)}{(1+|\tebu|)^{1+\gamma}}\, dx + \frac{1}{\gamma}\right).
\end{equation} 
It is worth to mention that $\bar c=\bar{c}(\data,\gamma,\delta)>0$ depends on $\data, \gamma$ and $\delta$ and it is an increasing function of $\gamma$ and $\delta$.
Since $\gamma\in (0,1/(s_G n))$ is fixed, we may choose $\alpha$ such that \begin{equation}
    \label{alpha-range} \alpha \in \left(1,\frac{n}{n-\kappa}\right) \qquad \text{and} \qquad  \frac{\alpha s_G\gamma}{\alpha-1} \leq 1.
\end{equation} Then
\[\int_{B_{1}}(1+|\tebu|)^{\frac{\alpha s_G\gamma}{\alpha-1}} \, dx\leq \int_{B_{1}} 1+|\tebu|\,dx\leq 1+|B_1|\] and we can estimate 
\begin{flalign}\label{int-vt-2} &\int_{B_{1}}\frac{\bG\left((1+|\tebu)|)^{1+\gamma}|D\phi|\right)}{(1+|\tebu|)^{1+\gamma}}\, dx \leq \int_{B_{1}}\frac{\bG\left((1+|\tebu|)|D\phi|\right)}{(1+|\tebu|)^{1+\gamma-\gamma s_G}}\, dx \notag \\
&\qquad\qquad\qquad\leq\left(\int_{B_{1}}(1+|\tebu|)^{\frac{\alpha s_G\gamma}{\alpha-1}} \, dx\right)^\frac{\alpha-1}{\alpha}\left(\int_{B_{1}} \left(\frac{\bG\left((1+|\tebu|)|D\phi|\right) }{(1+|\tebu|)^{1+\gamma}}\right)^{\alpha} \, dx\right)^\frac{1}{\alpha}\notag \\
&\qquad\qquad\qquad\leq c\left(\int_{B_{1}} \left(\frac{\bG\left((1+|\tebu|)|D\phi|\right) }{(1+|\tebu|)^{1+\gamma}}\right)^{\alpha} \, dx\right)^\frac{1}{\alpha}.
\end{flalign}

Thus, from \eqref{int-vt} and \eqref{int-vt-2} we obtain 
\color{black}
\begin{equation}\label{do-int-2} \left(\int_{B_1}\left( \frac{\bG((1+|\tebu|)\phi^q)}{(1+|\tebu|)^{1+\gamma}}\right)^\frac{n}{n-\kappa}\,dx\right)^\frac{n-\kappa}{n}\leq c\left(\int_{B_{1}}\left(\frac{\bG\left((1+|\tebu|)|D\phi|\right) }{(1+|\tebu|)^{1+\gamma}}\right)^{\alpha}\, dx\right)^\frac{1}{\alpha} + c .
\end{equation}
For $7/8 \leq r_1  < r_2\leq 1$, we take a cut-off function $\phi$ satisfying \[\phi\equiv 1\quad\text{on }\ B_{r_1}\qquad\text{and}\qquad |D\phi|\leq\frac{100}{r_2-r_1}.\]
It then follows from the doubling property of $\bG$ and Lemma~\ref{lem:self} that for any $\upsilon \in (0,1/s_G)$
\begin{flalign*} \left(\int_{B_{7/8}}\left( \frac{\bG(1+|\tebu|)}{(1+|\tebu|)^{1+\gamma}}\right)^\frac{n}{n-\kappa}\,dx\right)^\frac{n-\kappa}{n}
& \leq c\left(\int_{B_{1}}\left(\frac{\bG\left(1+|\tebu|\right) }{(1+|\tebu|)^{1+\gamma}}\right)^{\upsilon}\, dx\right)^\frac{1}{\upsilon} + c \\
& \leq c \left(\int_{B_{1}} (\bG\left(1+|\tebu|\right) )^{\upsilon}\, dx\right)^\frac{1}{\upsilon} + c\\
& \leq \bG \left(\int_{B_{1}} (1+|\tebu|)\, dx \right) + c.
\end{flalign*}
In the last line, we have used Jensen's inequality with the concave function $t \mapsto \bG(t)^{\upsilon}$.
Recalling \eqref{bu-male}, we obtain 
\begin{flalign}\label{do-int-3}
\left(\int_{B_{7/8}}\left( \frac{\bG(1+|\tebu|)}{(1+|\tebu|)^{1+\gamma}}\right)^\frac{n}{n-\kappa}\,dx\right)^\frac{n-\kappa}{n} \leq c
\end{flalign}
for some $c=c(\data,\gamma)>0.$

To proceed further, we recall $H_s$ defined in~\eqref{Hs-def} with $s$ from~\eqref{s-range}. By~\eqref{Hs-prop-1} function $H_s$ satisfies also $H_s(Mt) \approx_{\data} \bg(t)^{1-s} \bG(t)^{s} $. We apply Young's inequality, \eqref{post-cava-2}, and \eqref{do-int-3} with suitable choice of $\phi$. In turn we see
\begin{flalign} \label{H-s}
\int_{B_{\rho_1}} & H_s(M|\teDbu|) \,dx
 \lesssim_{\data} \int_{B_{ \rho_1}} \frac{\bg(|\teDbu|)^{1-s}\bG(|\teDbu|)^{s}} {(1+|\tebu|)^{1+\gamma}}(1+|\tebu|)^{1+\gamma}\,dx \notag \\
& \leq \int_{B_{ \rho_1}} \frac{\bG(|\teDbu|)} {(1+|\tebu|)^{1+\gamma}}\,dx
+ \int_{B_{ \rho_1}} \frac{ \bg(|\teDbu|) (1+|\tebu|)^{\frac{1+\gamma} {1-s}}}{(1+|\tebu|)^{1+\gamma}} \,dx \notag \\
& \leq \int_{B_{ \rho_1}} \frac{\bG(|\teDbu|)} {(1+|\tebu|)^{1+\gamma}}\,dx
+ \int_{B_{ \rho_1}} \frac{ \bG(1+|\tebu|) (1+|\tebu|)^{\frac{(s + \gamma) s_G}{(1-s)}}}{(1+|\tebu|)^{1+\gamma}} \,dx \notag \\
&\notag \leq \int_{B_{ \rho_1}} \frac{\bG(|\teDbu|)} {(1+|\tebu|)^{1+\gamma}}\,dx\\
&\quad+ \left(\int_{B_{ \rho_1}}\left( \frac{\bG(1+|\tebu|)}{(1+|\tebu|)^{1+\gamma}}\right)^\frac{n}{n-\kappa}\,dx\right)^\frac{n-\kappa}{n}
\left( \int_{B_{ \rho_1}}  (1+|\tebu|)^\frac{(s + \gamma) s_G n}{(1-s)\kappa}\,dx \right)^\frac{\kappa}{n}
\end{flalign}
In order to use \eqref{bu-male} we need to have 
\begin{equation}\label{s1}
\frac{(s + \gamma) s_G n}{(1-s)\kappa} \leq 1. 
\end{equation}
Observe that by Remark \ref{rem:iG-male} we have
$2-i_G <({i_G-\gamma s_G n})/({i_G+ s_Gn}),$
and we can choose $\kappa \in [1, i_G)$ such that
\begin{equation*}
    \max\{2-i_{G},0\} < s \leq {\frac{\kappa-\gamma s_G n}{\kappa+ s_Gn}}<s_{\rm m}
\end{equation*}
and the bound \eqref{s1} follows. Then \eqref{H-s} combined with \eqref{post-cava-2}, \eqref{int-vt-2} and \eqref{do-int-3}  implies for $s<s_{\rm m}$ that
\begin{equation}
    \label{H-apriori}\int_{B_{3/4}}H_s(M|\teDbu|)\,dx\leq c_{\rm ap}=c_{\rm ap}(\data,\gamma).
\end{equation}

\textbf{Step 4. Contradiction argument. } We state the counter--assumption. 

Namely, we assume that there exists $\ve$ and a sequences of balls $\{B_{r_j}(x_j)\}$ and almost $\opA$-harmonic maps $\{\teu_j\}\subset W^{1,G}(B_{r_j}(x_j),\Rm)$ such that \begin{equation}
    \label{u-male-j}\barint_{B_{r_j}(x_j)}|\teu_j|\,dx\leq Mr_j,\quad M\geq 1,
\end{equation}
\begin{equation}
    \label{Du-male-j}\left|\barint_{B_{r_j}(x_j)} 
    \opA(x,\teDu_j):\teDvp\,dx\right|\leq\frac{2^{-j}}{r_j}\|\tevp\|_{L^\infty(B_{r_j}(x_j),\Rm)}
\end{equation} 
for all $\tevp\in W^{1,G}_0(B_{r_j}(x_j),\Rm)\cap L^\infty(B_{r_j}(x_j),\Rm),$  but such that
    \begin{equation}\label{Du-daleko-Dv}
    \barint_{B_{r_j/2}(x_j)}H_s(|\teDu_j-\teDv|)\,dx
    > \ve
\end{equation}
whenever $\tev\in W^{1,G}(B_{r/2}(x_0),\Rm)$ is an $\opA$-harmonic map in  $ B_{r/2}(x_0)$ satisfying \begin{equation}\label{v-male-j}
    \barint_{B_{r_j/2}(x_0)}|\tev|\,dx\leq 2^nMr_j
    \quad\text{and}\quad 
    \left(\barint_{B_{r_j/2}(x_j)}H_s(M|\teDv|)\,dx\right)\leq c,
\end{equation} 
where $c=c(\data)>0.$ 

Let $\tebu$ be scaled as in~\eqref{bar}, but with the use of $x_j$ and $r_j,$ that is we set\[\tebu_j(x):=\frac{\teu(x_j+r_jx)}{Mr_j}\quad \text{and}\quad  \bopA(x_j+r_jx,\texi):=\opA(x_j+r_jx,M\texi).\]
In such a case by~\eqref{u-male-j} we get that\begin{equation}
    \label{buj-male} \barint_{B_1}|\tebu_j|\,dx\leq 1,
\end{equation}
so by~\eqref{Du-male-j} we infer that for all $\teeta=\tevp(x_j+r_jx)/r_j\in W^{1,G}_0(B_1,\Rm)\cap L^\infty(B_1,\Rm) $ it holds \begin{equation}
    \label{Dbu-male-j} \left|\barint_{B_1} 
    \bopA(x,\teDbu_j):\teDet\,dx\right|\leq {2^{-j}}\|\teeta\|_{L^\infty(B_1,\Rm)}
\end{equation}
and\begin{equation}
    \label{HsM}
    \barint_{B_1}H_s(M |\teDbu_j-\teDbv|)\,dx  > \ve
\end{equation}
whenever $\tebv\in W^{1,G}(B_{1/2}(x_0),\Rm)$ is an $\bopA$-harmonic map in  $ B_{1/2}(x_0)$ satisfying \begin{equation}\label{bv-male}
    \barint_{B_{1/2} }|\tebv|\,dx\leq 2^n 
    \quad\text{and}\quad 
    \left(\barint_{B_{1/2} }H_s(M|\teDbv|)\,dx\right)\leq c,
\end{equation} 
where $c=c(\data)>0.$

Since $\tebu$ satisfies~\eqref{buj-male}, we have~\eqref{H-apriori} for $s_{\rm o}<s_{\rm m}$  from~\eqref{s-range}. Therefore 
\[\int_{B_{3/4}}H_{s_{\rm o}}(M|\teDbu_j|)\,dx\leq C\quad\text{for}\quad C=C(\data,\gamma).\]
We fix any $s<s_{\rm o}$ from the range~\eqref{s-range}. Then we pick $\epsilon>0$ for which there exist $ c,t_0>0$, such that $H_s(t)\geq c t^{1+\epsilon}$ and $H_s(t)\geq c g^{1+\epsilon}(t)$ for all $t\geq t_0.$ In turn, we conclude with the following estimates uniform in $j$ \begin{equation}
    \label{unif-int-Duj}
    \int_{B_{3/4}}g^{1+\epsilon}(M|{\teDbu_j}|)\,dx\leq c_1\quad\text{and}\quad 
    \int_{B_{3/4}}(M|{\teDbu_j}|)^{1+\epsilon}\,dx\leq c_2
\end{equation}
with $c_1,c_2$ depending on $\data$ and $\gamma$ only. Further we infer that there exist \[\text{$\tewtu\in W^{1,H_s}(B_{3/4},\Rm),\quad$ $ \tfA\in L^{1+\epsilon}(B_{3/4},\rnm),\quad$ and $\quad\mathfrak{h}\in L^{H_s}(B_{3/4})$}\] such that up to a subsequence \begin{equation}
    \begin{split}\label{lots-of-conv}
        &{\teDbu_j}- {\teDwtu}\rightharpoonup 0\qquad \text{in } \ L^{H_s}(B_{3/4},\rnm),\\
        &|{\teDbu_j}-{\teDwtu}| \rightharpoonup \mathfrak{h}\qquad \text{in } \ L^{H_s}(B_{3/4}),\\
        &\bopA(x,{\teDbu_j}) \rightharpoonup \tfA\qquad \text{in } \ L^{1+\epsilon}(B_{3/4},\rnm),\\
        &\tebu_j\to \tewtu\quad \text{ strongly in  $\ L^{H_s}(B_{3/4},\Rm)\ $ and a.e. in $B_{3/4}$.} 
    \end{split}
\end{equation}
By~\eqref{buj-male},  lower semicontinuity of a functional $\tevp\mapsto\barint_{B_{1/2}}H_s(M|\teDvp|
)\,dx$, and~\eqref{H-apriori} we have \begin{equation}\label{wtu-male}
    \barint_{B_{1/2} }|\tewtu|\,dx\leq 2^n 
    \quad\text{and}\quad 
     \barint_{B_{1/2} }H_s(M|{\teDwtu}|)\,dx\leq c,
\end{equation}

\textbf{Step 5. Strong convergence of gradients. } Our aim is now to prove that  \begin{equation}
    \label{strong-conv-grad} {\teDbu_j}\to{\teDwtu}\qquad\text{in }\ L^{H_s}(B_{3/4},\rnm).
\end{equation}
For this we need to show that $\mathfrak{h}\in L^{H_s}(B_{3/4})$ from~\eqref{lots-of-conv} satisfies
\begin{equation}
    \label{h=0} \mathfrak{h}=0
\end{equation}
a.e. in  $B_{3/4}$.
This almost everywhere and weak convergence in $L^1$ implies strong $L^1$-convergence of ${\teDbu_j}\to{\teDwtu}$ in $B_{3/4}$.
Using the monotonicity property of $H_s$ and~\eqref{bv-male}, for sufficiently small $\wt\epsilon \in (0, \frac{1}{s_{H_{s}}})$ we have
\begin{align}\label{aux-esty1}
& \int_{B_{3/4}} H_s ( M |{\teDbu_j} - {\teDwtu}|) \,dx \\
&\quad \leq \int_{B_{3/4}} H_s( M |{\teDbu_j} - {\teDwtu}|)^{\wt\epsilon^2} H_s( M |{\teDbu_j}| + M |{\teDwtu}|)^{1-\wt\epsilon^2} \,dx \notag \\
&\quad \leq \bigg( \int_{B_{3/4}} H_s( M |{\teDbu_j} - {\teDwtu}|)^{\wt \epsilon} \,dx \bigg)^{\wt\epsilon} \bigg( \int_{B_{3/4}} H_s( M |{\teDbu_j}| + M |{\teDwtu}|)^{1 + \wt\epsilon} \,dx \bigg)^{1-\wt\epsilon} \notag \\
&\quad \leq c \bigg( \int_{B_{3/4}} H_s( M |{\teDbu_j} - {\teDwtu}|)^{\wt\epsilon} \,dx \bigg)^{\wt\epsilon}. \notag
\end{align}
Denoting
$$\Psi (t) = \int_{0}^{t} \frac{H_{s}^{-1}(\tau^{1/\wt\epsilon})}{\tau} \, d \tau,$$
one can immediately check
$$\frac{t \Psi''(t)}{\Psi'(t)} 
= \frac{t^{1/\wt \epsilon}}{\wt \epsilon H_{s}' ( H_{s}^{-1}(t^{1/\wt \epsilon}) ) H_{s}^{-1}(t^{1/\wt \epsilon}) } - 1
\geq \frac{1}{\wt \epsilon s_{H_{s}}} -1 >0,$$
and so $\Psi$ is a Young function.
We then apply Jensen's inequality to \eqref{aux-esty1} to obtain
$$\int_{B_{3/4}} H_s ( M |{\teDbu_j} - {\teDwtu}|) \,dx 
\leq c \bigg[ H_{s} \bigg( M  \int_{B_{3/4}} |{\teDbu_j} - {\teDwtu}| \,dx \bigg) \bigg]^{\wt\epsilon^2} \stackrel{j \to \infty}{\longrightarrow} 0.$$
Hence, it remains to show \eqref{h=0} to obtain \eqref{strong-conv-grad}.

We pick $\bar{x}$ being a Lebesgue's point simultaneously for ${\tewtu}, {\teDwtu}, \mathfrak{h},\tfA$, that is
\begin{flalign}
    \label{Leb-point}
    \lim_{\vr\to 0}\barint_{B_\vr(\bar{x})}&H_s(M|{\tewtu}-{\tewtu}(\bar{x})|)+H_s(M|{\teDwtu}-{\teDwtu}(\bar{x})|)\\
    &+H_s(M|\mathfrak{h}-\mathfrak{h}(\bar{x})|)+|\tfA-\tfA(\bar{x})|^{1+\epsilon}\,dx=0\nonumber
\end{flalign}
and
\begin{equation}
    \label{fin-val-barx}
    |{\tewtu}(\bar{x})|+|{\teDwtu}(\bar{x})|+|\mathfrak{h}(\bar{x})|+|\tfA(\bar{x})|<\infty.
\end{equation}
Almost every point of $B_{3/4}$ satisfies this conditions. Thus it is enough to  show that~\eqref{h=0} holds for $\bar{x}$.

We restrict our attention to $\vr$ small enough for $B_\vr(\bar x)\subset B_{3/4}$ and we set the linearization of ${\tewtu}$ at $\bar{x}$\begin{equation}
    \label{l-vr} \teelvr(x):=({\tewtu})_{B_{\vr}(\bar x)} + {\teDwtu}(\bar{x}):(x-\bar{x}).
\end{equation}Having the classical Poincar\'e inequality and~\eqref{Leb-point}, we obtain that
\begin{equation}
    \label{poinc-zero}
\lim_{\vr\to 0}\barint_{B_\vr(\bar x)}\left|\frac{{\tewtu}-\teelvr}{\vr}\right|^{1+\epsilon}\,dx\leq c \lim_{\vr\to 0}\barint_{B_\vr(\bar x)}\left|{\teDwtu}-{\teDwtu}(\bar{x})\right|^{1+\epsilon}\,dx=0.
\end{equation} 

Let us set\[\tI_{j,\vr}^0:=\barint_{B_{\vr/2}(\bar{x})}|{\teDbu_j}-{\teDwtu}|\,dx.\]
By~\eqref{lots-of-conv} we have the weak convergence of $|{\teDbu_j}-{\teDwtu}|\to\mathfrak{h}$ in $L^{H_s}(B_\vr).$  Since $\bar{x}$ is a~Lebesgue's point of $\teDwtu$ we infer that\begin{equation}
    \label{h-of-barx}\mathfrak{h}(\bar{x})=\lim_{\vr\to 0}\lim_{j\to\infty}\tI_{j,\vr}^0.
\end{equation}
In order to prove that $\mathfrak{h}(\bar{x})=0$, let us write\begin{flalign}
    \nonumber \tI_{j,\vr}^0=\barint_{B_{\vr/2}(\bar{x})}|{\teDbu_j}-{\teDwtu}|\,dx&=\barint_{B_{\vr/2}(\bar{x})}\mathds{1}_{\{|\tebu_j-\teelvr|\geq\vr\}}|{\teDbu_j}-{\teDwtu}|\,dx\\
    &\ \ + \barint_{B_{\vr/2}(\bar{x})}\mathds{1}_{\{|\tebu_j-\teelvr|<\vr\}}|{\teDbu_j}-{\teDwtu}|\,dx:=\tI^1_{j,\vr}+\tI^2_{j,\vr}\label{2ndline}
\end{flalign}
and prove the convergence of both terms first when $j\to \infty$ and then $\vr\to 0$.

We start with $\tI^1_{j,\vr}$. Let us observe that\begin{flalign*}
    \tI^1_{j,\vr}&\leq\barint_{B_{\vr/2}(\bar{x})}\mathds{1}_{\{|\tebu_j-{\tewtu}|\geq\vr\}}|{\teDbu_j}-{\teDwtu}|\,dx\\
    & \quad +\barint_{B_{\vr/2}(\bar{x})}\mathds{1}_{\{|\tewtu-\teelvr|\geq\vr\}}|{\teDbu_j}-{\teDwtu}|\,dx=:\tI^{1,1}_{j,\vr}+\tI^{1,2}_{j,\vr}. 
\end{flalign*}
Notice that  $\tI^{1,1}_{j,\vr}$ vanishes as $j\to\infty$. Indeed, $\tI^{1,1}_{j,\vr}\geq 0$ and by H\"older inequality we have\begin{flalign*}
    \tI^{1,1}_{j,\vr}&\leq \frac{1}{|B_{\vr/2}(\bar{x})|} \left(\int_{B_{\vr/2}(\bar{x})}|{\teDbu_j}-{\teDwtu}|^{1+\epsilon}\,dx\right)^{\frac{1}{1+\epsilon}}\left(|\{x\in B_{3/4}:\ |\tebu_j-{\tewtu}|\geq\vr/2\}|\right)^\frac{\epsilon}{1+\epsilon}. 
\end{flalign*} 
Since by~\eqref{lots-of-conv} one has that $|\tebu_j-{\tewtu}|\to 0$ strongly in $L^1(B_{3/4})$, so \[\lim_{j\to\infty} |\{x\in B_{3/4}:\ |\tebu_j-{\tewtu}|\geq\vr/2\}|=0.\] The rest of the terms are bounded as $\epsilon$ is chosen such that~\eqref{unif-int-Duj} is true and $|{\teDwtu}|$ shares the same a priori estimates as $|\teDbu_j|$. Therefore, we infer that $\lim_{j\to\infty} \tI^{1,1}_{j,\vr}=0$. On the other hand, $\tI^{1,2}_{j,\vr}$ is convergent when $j\to\infty$, because of the weak convergence of $|{\teDbu_j}-{\teDwtu}|\to\mathfrak{h}$ in $L^{H_s}(B_\vr)$. Hence, we get
\[\lim_{j\to\infty}\tI^{1,2}_{j,\vr}=\barint_{B_{\vr/2}(\bar{x})}\mathds{1}_{\{|\tewtu-\teelvr|\geq\vr\}}\mathfrak{h}\,dx=:\tI^{1,2}_\vr.\]
We can estimate further
\begin{flalign*}
    \tI^{1,2}_\vr&\leq \left(\barint_{B_\vr(\bar x)}\mathfrak{h}^{1+\epsilon}\,dx\right)^\frac{1}{1+\epsilon}\left(\barint_{B_\vr(\bar x)}\mathds{1}_{\{|\tewtu-\teelvr|\geq\vr/2\}}\,dx\right)^\frac{\epsilon}{1+\epsilon}\\
    &\leq c\left[\left(\barint_{B_\vr(\bar x)}|\mathfrak{h}-\mathfrak{h}(\bar{x})|^{1+\epsilon}\,dx\right)^\frac{1}{1+\epsilon}+\mathfrak{h}(\bar{x})\right]\left(\barint_{B_\vr(\bar x)}\left|\frac{\tewtu-\teelvr}{\vr}\right|^{1+\epsilon}\,dx\right)^\frac{\epsilon}{1+\epsilon}
\end{flalign*}
which tends to $0$ as $\vr\to 0$ as $\bar{x}$ is a Lebesgue's point of $\mathfrak{h}$ as in~\eqref{Leb-point} and the last bracket converges to $0$ due to~\eqref{poinc-zero}. Altogether, we have that  $\tI^{1}_{j,\vr}$ vanishes in the limit, so we will now concentrate on $\tI^2_{j,\vr}$ for which we have
\begin{flalign*}
    \tI^2_{j,\vr}&\leq \barint_{B_{\vr/2}(\bar{x})}\mathds{1}_{\{|\tebu_j-\teelvr|<\vr\}}|{\teDbu_j}-D{\teelvr}|\,dx+2^n\barint_{B_{\vr}(\bar{x})}\mathds{1}_{\{|\tebu_j-\teelvr|<\vr\}}|D{\teelvr}-{\teDwtu}|\,dx\\
    &=:\tI^{2,1}_{j,\vr}+\tI^{2,2}_{j,\vr}.
\end{flalign*}
By~\eqref{lots-of-conv},~\eqref{Leb-point}, and~\eqref{l-vr} we have that $\lim_{\vr\to 0}\limsup_{j\to\infty}\tI^{2,2}_{j,\vr}=0.$ Proving the convergence \begin{equation}
    \label{I22-conv}\limsup_{\vr\to 0}\limsup_{j\to\infty}\tI^{2,1}_{j,\vr}=0
\end{equation}
 requires more arguments. We take \[\text{$\phi\in C_c^\infty(B_\vr(\bar x))\quad$ with $0\leq \phi\leq 1$, $\ \phi \equiv 1$ on $B_{\vr/2}(\bar x)\ $ and $\ |D\phi|\leq 4/\vr.$}\] Let \[\teeta=\phi T_\vr ({\tebu_j}-\teelvr),\] 
where the truncation is defined in~\eqref{Tk}. Let us denote 
\begin{equation}
    \label{PjP}\Pj:=\frac{({\tebu_j}-\teelvr)\otimes ({\tebu_j}-\teelvr)}{|{\tebu_j}-\teelvr|^2}\qquad\text{and}\qquad \P:=\frac{({\tewtu}-\teelvr)\otimes ({\tewtu}-\teelvr)}{|{\tewtu}-\teelvr|^2}
\end{equation}
when $|{\tebu_j}-\teelvr|\neq 0$ and $|{\tewtu}-\teelvr|\neq 0,$ respectively. 
Within this notation we have that\begin{flalign*}
    &\left(\bopA(x,{\teDbu_j})-\bopA(x,D{\teelvr})\right):\teDet\\
    &\qquad=\mathds{1}_{\{|{\tebu_j}-\teelvr|<\vr\}}\left[\left(\bopA(x,{\teDbu_j})-\bopA(x,D{\teelvr})\right):{D({\tebu_j}-\teelvr)}\right]\phi\\
    &\qquad\quad+\frac{\vr\mathds{1}_{\{|{\tebu_j}-\teelvr|\geq\vr\}}}{|{\tebu_j}-\teelvr|}\left[\left(\bopA(x,{\teDbu_j})-\bopA(x,D{\teelvr})\right):(\id-\Pj){D({\tebu_j}-\teelvr)}\right]\phi\\
    &\qquad\quad+\left(\bopA(x,{\teDbu_j})-\bopA(x,D{\teelvr})\right):\left[T_\vr({\tebu_j}-\teelvr)\otimes D\phi\right] \\
    &\qquad=: G_{j,\vr}^{1}(x)+G_{j,\vr}^2(x)+G_{j,\vr}^3(x) .
\end{flalign*}
Moreover,  we recall that  since $\teelvr$ is affine, whenever $B\subset B_1$ it holds that\[\int_B\bopA(x,D{\teelvr}):\teDvp\,dx=0\qquad\text{for every }\ \tevp\in W^{1,1}_0(B,\Rm).\]
Therefore, by~\eqref{Dbu-male-j} it is justified to write that\begin{equation}
    \label{4.46}
    0\leq\barint_{B_\vr(\bar x)}G^1_{j,\vr}(x)\,dx\leq 2^{-j}\vr^{1-n}-\barint_{B_\vr(\bar x)}G^2_{j,\vr}(x)\,dx-\barint_{B_\vr(\bar x)}G^3_{j,\vr}(x)\,dx.
\end{equation}The first term in the above display is nonnegative because of the monotonicity of~$\bopA.$ Instrumental for proving that $\tI_{j,\vr}^{2,1}\to 0$ is to establish that \begin{equation}
    \label{G1-conv} \limsup_{\vr\to 0}\limsup_{j\to\infty}\barint_{B_\vr(\bar x)}G^1_{j,\vr}(x)\,dx=0,
\end{equation}
which will be proven provided  one justifies that the last two terms of~\eqref{4.46} vanish in the limit. We will show first that\begin{equation}\label{G2-conv}
\limsup_{\vr\to 0}    \limsup_{j\to\infty} \left(-\,\barint_{B_\vr(\bar x)} G^2_{j,\vr}\,dx\right)\leq 0.
\end{equation}
The quasi-diagonal structure of $\bopA$ ensures that $\bopA(x,{\teDbu_j}):[(\id-\Pj){\teDbu_j}]\geq 0$, see~\eqref{calc-eta}. Therefore,\begin{flalign}
    \nonumber &\big(\bopA(x,{\teDbu_j})-\bopA(x,D{\teelvr})\big):(\id-\Pj){D({\tebu_j}-\teelvr)}\\
    &\qquad\geq -\bopA(x,{\teDbu_j}):(\id-\Pj)D{\teelvr}-\bopA(x,D{\teelvr}):(\id-\Pj){D({\tebu_j}-\teelvr)}.\label{monot-ell-vr}
\end{flalign}
Recall that $\Pj$ and $\P$, defined in~\eqref{PjP}, are bounded. Notice that for $j\to\infty$ we have $\mathds{1}_{\{|{\tebu_j}-\teelvr|\geq\vr\}}\Pj\to \mathds{1}_{\{|{\tewtu}-\teelvr|\geq\vr\}}\P $ almost everywhere and thus, by the Lebesgue's dominated convergence theorem,  also strongly in $L^t(B_{3/4})$ for every $t\geq 1.$ Moreover,  $\mathds{1}_{\{|{\tebu_j}-\teelvr|\geq\vr\}}|{\tebu_j}-\teelvr|^{-1}\to \mathds{1}_{\{|{\tewtu}-\teelvr|\geq\vr\}}|{\tewtu}-\teelvr|^{-1}$ almost everywhere and, as a uniformly bounded sequence of functions, it converges also strongly in $L^t(B_{3/4})$ for every $t\geq 1.$ Having this,~\eqref{monot-ell-vr}, and~\eqref{lots-of-conv}, we obtain\begin{flalign}
    \nonumber \limsup_{j\to\infty} &\left(-\barint_{B_\vr(\bar x)} G^2_{j,\vr}\,dx\right)\leq \barint_{B_\vr(\bar x)}\tfA:(\id-\P)D{\teelvr}\frac{\vr\mathds{1}_{\{|{\tewtu}-\teelvr|\geq\vr\}}}{|{\tewtu}-\teelvr|}\,dx\\
    &\qquad\qquad+\barint_{B_\vr(\bar x)}\bopA(x,D{\teelvr}):(\id-\P){D({\tewtu}-\teelvr)}\frac{\vr\mathds{1}_{\{|{\tewtu}-\teelvr|\geq\vr\}}}{|{\tewtu}-\teelvr|}\,dx=:\tII^1_\vr+\tII^2_\vr.
\end{flalign}
We can estimate
\begin{flalign*}
    |\tII^1_{\vr}|&\leq c\,\barint_{B_\vr(\bar{x})}|\tfA|\left|\frac{{\tewtu}-\teelvr}{\vr}\right|^{\epsilon}\,dx\leq c\left(\barint_{B_\vr(\bar{x})}|\tfA|^{1+\epsilon}\,dx\right)^{\frac{1}{1+\epsilon}}\left(\barint_{B_\vr(\bar{x})}\left|\frac{{\tewtu}-\teelvr}{\vr}\right|^{1+\epsilon}\,dx\right)^{\frac{\epsilon}{1+\epsilon}}
\end{flalign*}
where the first term is bounded and the second convergent to zero by~\eqref{poinc-zero}.
On the other hand, by~\eqref{l-vr} and \eqref{fin-val-barx} we may estimate
\begin{flalign*}
    |\tII^2_{\vr}|&\leq c\,\barint_{B_\vr(\bar{x})}|{D({\tewtu}-\teelvr)}|\left|\frac{{\tewtu}-\teelvr}{\vr}\right|^\epsilon\,dx\\
    &\leq c\left(\barint_{B_\vr(\bar{x})}|{\teDwtu}-{\teDwtu}(\bar{x})|^{1+\epsilon}\,dx\right)^\frac{1}{1+\epsilon}\left(\barint_{B_\vr(\bar{x})}\left|\frac{{\tewtu}-\teelvr}{\vr}\right|^{1+\epsilon}\,dx\right)^\frac{\epsilon}{1+\epsilon}
\end{flalign*}
where, again, the first term is bounded and the second convergent to zero by~\eqref{poinc-zero}. Summing up the information from the last three displays we get~\eqref{G2-conv}.

Now we concentrate on justifying that\begin{equation}
    \label{G3-conv} \lim_{\vr\to 0}\lim_{j\to\infty}\left| \barint_{B_\vr(\bar x)}G^3_{j,\vr}(x)\,dx\right|=0.
\end{equation} Let us observe that because~\eqref{lots-of-conv} provides weak convergence of $\bopA(x,{\teDbu_j})$ to $\tfA$ in $L^{1+\epsilon}$ and strong convergence ${\tebu_j}$ to ${\tewtu}$  in $L^{1+\epsilon}$, we obtain that
\[\lim_{j\to\infty} \barint_{B_\vr(\bar x)}G^3_{j,\vr}(x)\,dx=\barint_{B_\vr(\bar x)}\big(\tfA-\bopA(x,D{\teelvr})\big):[T_\vr({\tewtu}-\teelvr)\otimes D\phi]\,dx.\]
By H\"older inequality and the choice of $\phi$, we estimate further \begin{flalign*}
    \lim_{j\to\infty} \barint_{B_\vr(\bar x)}G^3_{j,\vr}(x)\,dx &\leq  \left(\barint_{B_\vr(\bar x)}|\tfA-\tfA(\bar{x})|^{1+\epsilon}+|\tfA(\bar{x})|^{1+\epsilon} + g(|{\teDwtu}(\bar{x})|)^{1+\epsilon} \, dx \right)^\frac{1}{1+\epsilon}\\
    &\quad \cdot \left(\barint_{B_\vr(\bar x)}\left(\frac{\min\{\vr,|{\tewtu} -\teelvr|\}}{\vr}\right)^{\frac{1+\epsilon}{\epsilon}}\,dx\right)^{\frac{\epsilon}{1+\epsilon}},
\end{flalign*}
where the first integral on the right-hand side is finite and the second term converges to zero. Indeed, since $0 \leq \min\{ \rho, |{\tewtu} - l_{\rho}|\} \leq \rho$, we have
\begin{flalign*} \barint_{B_\vr(\bar x)}\left(\frac{\min\{\vr,|{\tewtu} -\teelvr|\}}{\vr}\right)^{\frac{1+\epsilon}{\epsilon}}\,dx & \leq  \barint_{B_\vr(\bar x)}\left(\frac{\min\{\vr,|{\tewtu} -\teelvr|\}}{\vr}\right)^{1+\epsilon}\,dx \\
&\leq  \barint_{B_\vr(\bar x)}\left(\frac{|{\tewtu} -\teelvr|}{\vr}\right)^{1+\epsilon}\,dx  \xrightarrow[\vr\to 0]{}0,\end{flalign*}
where the last convergence results from~\eqref{poinc-zero}. Therefore, we get~\eqref{G3-conv}.

We have shown~\eqref{G2-conv} and \eqref{G3-conv}, so because of~\eqref{4.46} the limit  \eqref{G1-conv} follows. Hence, we are in the position to prove~\eqref{I22-conv}. In the view of~\eqref{opA:strict-monotonicity}, \eqref{G1-conv} implies that
\begin{equation}\label{pre-grad-conv} 
    \limsup_{\vr\to\infty}\limsup_{j\to\infty}\barint_{B_\vr(\bar x)}\mathds{1}_{\{|{\tebu_j}-\teelvr|<\vr\}}  \frac{g(|{\teDbu_j}|+|D{\teelvr}|)}{|{\teDbu_j}|+|D{\teelvr}|}|{D({\tebu_j}-\teelvr)}|^2\,dx=0.
\end{equation}
At this stage, we calculate similarly to \eqref{aux-esty1} in order to show \eqref{I22-conv}.
For any $\hat\epsilon<\frac{2}{s_{G}}$, it is readily checked that $t \mapsto t g(t)^{-\hat\epsilon/(2- \hat\epsilon)}$ is a monotone increasing function. 
Then
\begin{flalign}
   & \nonumber \barint_{B_\vr(\bar x)}\mathds{1}_{\{|{\tebu_j}-\teelvr|<\vr\}}|{D({\tebu_j}-\teelvr)}|\,dx\\
   & \nonumber \leq \barint_{B_\vr(\bar x)}\mathds{1}_{\{|{\tebu_j}-\teelvr|<\vr\}} \left(|{D({\tebu_j}-\teelvr)}|^2\frac{g(|{\teDbu_j}|+|D{\teelvr}|)}{|{\teDbu_j}|+|D{\teelvr}|}\right)^{\frac{\hat\epsilon}{2}} \frac{(|{\teDbu_j}|+|D{\teelvr}|)^{1-\frac{\hat\epsilon}{2}}}{{g}(|{\teDbu_j}|+|D{\teelvr}|)^{\frac{\hat\epsilon}{2}}}\,dx\\
   & \nonumber \leq \left(\barint_{B_\vr(\bar x)}\mathds{1}_{\{|{\tebu_j}-\teelvr|<\vr\}} |{D({\tebu_j}-\teelvr)}|^2\frac{g(|{\teDbu_j}|+|D{\teelvr}|)}{|{\teDbu_j}|+|D{\teelvr}|} \,dx \right)^{\frac{\hat\epsilon}{2}}\\
   & \nonumber \quad \cdot \left( \barint_{B_\vr(\bar x)}\mathds{1}_{\{|{\tebu_j}-\teelvr|<\vr\}} \frac{|{\teDbu_j}|+|D{\teelvr}|}{g(|{\teDbu_j}|+|D{\teelvr}|)^{\frac{\hat\epsilon}{2-\hat\epsilon}}}\,dx \right)^{1-\frac{\hat\epsilon}{2}}\\
   &\nonumber \leq c \left(\barint_{B_\vr(\bar x)}\mathds{1}_{\{|{\tebu_j}-\teelvr|<\vr\}} |{D({\tebu_j}-\teelvr)}|^2 \frac{g(|{\teDbu_j}|+|D{\teelvr}|)}{|{\teDbu_j}|+|D{\teelvr}|} \,dx\right)^\frac{\hat\epsilon}{2}\\ 
   &\nonumber \quad\cdot
   \left(\barint_{B_\vr(\bar x)} (|{\teDbu_j}|+|D{\teelvr}| + 1) \,dx\right)^{1-\frac{\hat\epsilon}{2}},
\end{flalign}
where $c=c(g)>0$. 
Noting that the very last term in the above display is bounded, by~\eqref{pre-grad-conv} we infer that~\eqref{I22-conv} holds.

Summing up all the convergences of this step, we get in~\eqref{h-of-barx} that $\mathfrak{h}(\bar{x})=0$ and, consequently,~\eqref{h=0} holds almost everywhere in $B_{3/4}.$ As explained in the beginning of this step, this suffices to get the final aim of \textbf{Step~5}, that is strong convergence of gradients~\eqref{strong-conv-grad}.

\medskip

\textbf{Step 6. $\bopA$-harmonicity of the limit map and conclusion by contradiction.} Having~\eqref{strong-conv-grad}, we can pass to the limit in~\eqref{Dbu-male-j} with $j\to\infty$ getting that\begin{equation}
    \label{wt-u-harm} \barint_{B_{1/2}} \bopA(x,{\teDwtu}):\teDet\,dx=0\qquad\text{
for $\ \ \teeta\in C_c^\infty(B_{1/2},\Rm)$.}
\end{equation} Therefore, if ${\tewtu}\in W^{1,G}(B_{1/2},\Rm)$, then it will be proven to be  $\bopA$-harmonic. Indeed, since we know~\eqref{wtu-male}, it is allowed to take $\tebv={\tewtu}$ in \eqref{HsM}. Note that in such a case~\eqref{wtu-male} is precisely the restriction on the test function from~\eqref{bv-male}. Then, in the view of~\eqref{strong-conv-grad}, taking $j$ large enough, we will get the desired contradiction. Hence, it remains to prove that $|{\teDwtu}|\in L^G(B_{1/2})$.

We have~\eqref{do-int-3} for each ${\tebu_j}$ with the constant independent of $j$, so by the lower semicontinuity we can write that \begin{flalign*}
\left(\int_{B_{7/8}}\left( \frac{\bG(1+|{\tewtu}|)}{(1+|{\tewtu}|)^{1+\wtgamma}}\right)^\frac{n}{n-\kappa}\,dx\right)^\frac{n-\kappa}{n} &\leq\liminf_{j\to\infty}\left(\int_{B_{7/8}}\left( \frac{\bG(1+|{\tebu_j}|)}{(1+|{\tebu_j}|)^{1+\wtgamma}}\right)^\frac{n}{n-\kappa}\,dx\right)^\frac{n-\kappa}{n}\\
&\leq c
\end{flalign*}
for some $c=c(\data,\wtgamma)>0.$ Analogically, by~\eqref{st2-apriori} for ${\tebu_j}$ and with $\delta=2^{-j}$, by letting $j\to\infty$, Fatou's lemma on the left-hand side of the resultant inequality and \eqref{strong-conv-grad} on its right-hand side, we get \begin{flalign}\nonumber
  \int_{B_{3/4}\cap \{|{\tewtu}|< t\}}\bG(|{\teDwtu}|) \phi^q\,dx\leq &\,c_* \int_{B_{3/4}\cap \{|{\tewtu}|< t\}}\bG(|{\tewtu}|\,|D\phi|) \,dx\\
  &+\,c_*\,t \int_{B_{3/4}\cap \{|{\tewtu}|\geq t\}} {\bg(|{\teDwtu}|)}\phi^{q-1}| D\phi|\,dx\nonumber
\end{flalign} 
for every $t>0$ and $\phi\in C_c^\infty(B_{3/4})$ with $\phi\geq 0$, and $c_*=c_*(\data,q).$ We proceed as in the beginning of {\bf Step 3}. We multiply the above display by $(1+t)^{-(\wtgamma+1)},$ $\wtgamma>0$ to be chosen sufficiently small in a few lines, integrate it from zero to infinity and apply Cavalieri's principle (Lemma~\ref{lem:cava}) twice (with $\nu_1=\bG(|\teDwtu|) \phi^q$ and $\nu_2=\bG(|\tewtu|\,|D\phi|)$. Altogether we get
 \begin{flalign}\nonumber
 \frac{1}{\wtgamma} \int_{B_{3/4}}&\frac{\bG(|{\teDwtu}|) \phi^q}{(1+|{\tewtu}|)^{\wtgamma}}\,dx\leq
 \frac{c_*}{\wtgamma} \int_{B_{3/4}}\frac{\bG(|{\tewtu}|\,|D\phi|) }{(1+|{\tewtu}|)^{\wtgamma}}\,dx\\
  &+ c_* \int_0^\infty \frac{1}{(1+t)^{\wtgamma}}\int_{B_{3/4}\cap\{|{\tewtu}|\geq t\} } {\bg(|{\teDwtu}|)}\phi^{q-1}| D\phi|\,dx \,dt=\tIII_1+\tIII_2.\nonumber
\end{flalign}
To estimate further the very last term we note that $q$ is large enough to satisfy $s_G'\geq q'$, and so Lemma~\ref{lem:equivalences} implies that ${\bG}^*(\phi^{q-1}\bg(t))\leq c_G \phi^q \bG(t).$ Then,  using Young inequality~\eqref{in:Young} and by taking $\wtgamma \in (0,1/(2c_*c_G+1)],$ get
\begin{equation*}
\begin{split}\tIII_2&\leq \frac{c_*}{1-\wtgamma}\int_{B_{3/4} } {\bg(|{\teDwtu}|)}{(1+|{\tewtu}|)^{1-\wtgamma}}\phi^{q-1}| D\phi|\,dx\\
&\leq  \frac{1}{2\wtgamma c_G}  \int_{B_{3/4}}\frac{{\bG^*}\left( {\phi^{q-1}\bg(|{\teDwtu}|) }\right)}{(1+|{\tewtu}|)^{\wtgamma}} \, dx+ \frac{1}{2\wtgamma}\int_{B_{3/4}}\frac{\bG\left(\frac{2\wtgamma c_* c_G}{1-\wtgamma}(1+|{\tewtu}|)|D\phi|\right)}{(1+|{\tewtu}|)^{\wtgamma}}\,dx \\
& \leq \frac{1}{2\wtgamma} \int_{B_{3/4}}\frac{\bG\left(|{\teDwtu}| \right)\phi^{q }}{(1+|{\tewtu}|)^{\wtgamma}}  \,dx + \frac{c_* c_G}{1- \wtgamma} \int_{B_{3/4}}\frac{\bG\left((1+|{\tewtu}|) |D\phi|\right)}{(1+|{\tewtu}|)^{\wtgamma}}\, dx 
\end{split}
\end{equation*}
where in the last line we used that $\frac{2\wtgamma c_* c_G} {1-\wtgamma}<1$ can be taken out of the integrand by Jensen's inequality. Summing up we obtain 
 \begin{flalign}
 \int_{B_{3/4} }\frac{\bG(|{\teDwtu}|) \phi^q}{(1+|{\tewtu}|)^{\wtgamma}}\,dx&\leq
 \,{C} \int_{B_{3/4}}\frac{\bG\left((1+|{\tewtu}|)|D\phi|\right)}{(1+|{\tewtu}|)^{\wtgamma}}\, dx \label{post-cava-final-step}
\end{flalign}
where $C=C(\data)
>0$. Then similar calculations to \eqref{Dvt}-\eqref{int-vt} yield that
\begin{equation*}
    \left(\int_{B_{3/4}}\left( \frac{\bG((1+|\tebu|)\phi^q)}{(1+|\tebu|)^{\wtgamma}}\right)^\frac{n}{n-1}\,dx\right)^\frac{n-1}{n}\leq
 c \int_{B_{3/4}}\frac{\bG\left((1+|\tebu|)|D\phi|\right)}{(1+|\tebu|)^{\wtgamma}}\, dx.
\end{equation*}
holds with $c=c(\data)>0$. Indeed, in \textbf{Step 4}, we have checked that the above $c=c(\data,\wtgamma)$ is an increasing function of $\wtgamma$. As we consider small $\wtgamma$, $c$ in fact depends only on $\data$.
For $5/8 \leq r_{1} < r_{2} \leq 3/4$ we take $\phi \in C_{c}^{\infty}(B_{r_2})$ to satisfy
$$\phi \equiv 1 \quad \text{on} \quad B_{r_1} \qquad \text{and} \qquad |{D \phi}| \leq \frac{100}{r_2 - r_1}.$$
Then the doubling property of $\bG$ and Lemma \ref{lem:self} gives
\begin{equation}\label{to-est}
    \left(\int_{B_{5/8}}\left( \frac{\bG(1+|\tebu|)}{(1+|\tebu|)^{\wtgamma}}\right)^\frac{n}{n-1}\,dx\right)^\frac{n-1}{n}\leq c \left( \int_{B_{3/4}} \left(\frac{\bG(1+|\tebu|)}{(1+|\tebu|)^{\wtgamma}} \right)^{\frac{1}{2s_G}}\, dx \right)^{2s_G}.
\end{equation}
We now restrict ourselves to $\wtgamma \in\big(0,\min\{1/(2c_*c_G+1), \frac{i_G}{2}\}\big)$ and define
$$\Psi_{\wtgamma}(t) = \int_{0}^{t} \frac{1}{\tau} \left( \frac{\bG(\tau)}{\tau^{\wtgamma}} \right)^{\frac{1}{2s_G}} \, d \tau,$$
which is an increasing concave function on $[0,\infty)$ satisfying
$$\frac{t \Psi_{\wtgamma}''(t)}{\Psi_{\wtgamma}'(t)}=\frac{t g(t)}{2s_G \bG(t)} - 1 - \frac{\wtgamma}{2s_G} \in (-1,-1/2) \quad \text{and} \quad \Psi_{\wtgamma}(t) \approx \left(\frac{\bG(t)}{t^{\wtgamma}}\right)^{\frac{1}{2s_G}}.$$   Then Jensen's inequality gives 
\begin{flalign}
    \int_{B_{3/4}} &\left(\frac{\bG(1+|\tebu|)}{(1+|\tebu|)^{\wtgamma}} \right)^{\frac{1}{2s_G}}\, dx
    \leq c \int_{B_{3/4}} \Psi_{\wtgamma}(1+|\tebu|)\, dx \nonumber \leq c\, \Psi_{\wtgamma} \left( \int_{B_{3/4}} (1+|\tebu|) \, dx \right)\\
    &\qquad\quad \leq c\, \Psi_{\wtgamma} \left( \int_{B_{3/4}} (1+|\tebu|) \, dx + 1\right)
    \leq c\, \bG^\frac{1}{2s_G} \left( \int_{B_{3/4}} (1+|\tebu|) \, dx + 1\right)
    \leq c, \label{LG-final-step}
\end{flalign}
with $c=c(\data)>0.$ We used that if $t>1$ is arbitrarily fixed there exists $c>0$ independent of $\wtgamma$ such that for all $t>1$ and all $\wtgamma$, we have $\Psi_\wtgamma(t)\leq c \left({\bG(t)} \right)^{\frac{1}{2s_G}}$. 
By H\"older inequality,~\eqref{to-est} and \eqref{LG-final-step} we get that
\begin{equation}\label{to-est-2}
   \int_{B_{5/8}}\frac{\bG(1+|\tebu|)}{(1+|\tebu|)^{\wtgamma}}\,dx\leq \left(\int_{B_{5/8}}\left( \frac{\bG(1+|\tebu|)}{(1+|\tebu|)^{\wtgamma}}\right)^\frac{n}{n-1}\,dx\right)^\frac{n-1}{n}\leq c 
\end{equation}
with $c=c(\data)>0.$

We now consider \eqref{post-cava-final-step} with a cutoff function $\phi \in C_{c}^{\infty}(B_{5/8})$ satisfying
$$\phi \equiv 1 \quad \text{in} \quad B_{1/2} \qquad \text{and} \qquad |{D \phi}| \leq 100$$
and combine it with \eqref{to-est-2}, to obtain for some $c=c(\data)>0$ that
\begin{flalign*}
    \int_{B_{1/2} }\frac{\bG(|{\teDwtu}|)}{(1+|{\tewtu}|)^{\wtgamma}}\,dx&\leq
    {c} \int_{B_{3/4}}\frac{\bG(1+|{\tewtu}|)}{(1+|{\tewtu}|)^{\wtgamma}}\, dx \leq c.
\end{flalign*}
Therefore, using Fatou's lemma we justify that \begin{flalign*}
    \int_{B_{1/2} }{\bG(|{\teDwtu}|)}\,dx&\leq \limsup_{\wtgamma\to 0}
    \int_{B_{1/2} }\frac{\bG(|{\teDwtu}|)}{(1+|{\tewtu}|)^{\wtgamma}}\,dx\leq c.
\end{flalign*}
Consequently, we conclude that $|{D {\tewtu}}| \in L^{G}(B_{1/2})$.
This  completes the proof of Theorem \ref{theo:Ah-approx}.
\end{proof}

\section{Proof of Wolff potential estimates} \label{sec:mainproof}
\subsection{Comparison estimate} We need one more auxiliary estimate yielding comparison between energy of a weak solution and an $\opA$-harmonic function.
\begin{lem}\label{lem:A-h-appr} Under Assumption {\bf (A-vect)} suppose $u\in W^{1,G}(B_r,\Rm)$ is a weak solution to~\eqref{eq:mu} in $B_r=B_r(x_0),$ $r<1$ and  let $\ve\in(0,1)$. 
Then there exists a positive constant $c_{\rm s}=c_{\rm s}(\data,\ve)$ and a map $\tev$ being $\opA$-harmonic in $B_{r/2}$ and such that\begin{equation}
    \label{comp-est}\barint_{B_{r/2}}|\teDu-\teDv|\,dx\leq \frac{\ve}{r}\barint_{B_r}|\teu-(\teu)_{B_r}|\,dx+c_{\rm s}g^{-1}\left(\frac{|\temu|(B_r)}{r^{n-1}}\right).
\end{equation}  
\end{lem}
\begin{proof}Let us fix
\[\lambda:=\frac{1}{r}\barint_{B_r}|\teu-(\teu)_{B_r}|\,dx+g^{-1}\left(\delta \frac{|\temu|(B_r)}{r^{n-1}}\right)\]
with $\delta=\delta(\data,\ve)$ from Theorem~\ref{theo:Ah-approx} with $M=1$. If $\lambda=0,$ then $\teu$ is constant and $\tev=\teu$. Otherwise $\lambda>0$ and we can argue by scaling\[\tebu:=\frac{\teu-(\teu)_{B_r}}{\lambda},\qquad\bar{\temu}:=\frac{\temu}{g(\lambda)},\qquad\bopA(x,\texi)=\frac{\opA(x,\lambda\texi)}{g(\lambda)}.\]
Then
\[\left|\barint_{B_r} \bopA(x,\teDbu):\teDvp\,dx\right|\leq \frac{\|\tevp\|_{L^\infty(B_r)}|\temu|(B_r)}{g(\lambda) r^{n-1}}\leq\frac{\delta}{r}\|\tevp\|_{L^\infty(B_r)}.\]
By definition of $\tebu$ and $\lambda$ we notice that
\[\barint_{B_r}|\tebu|\,dx\leq r.\]
Therefore, by Theorem~\ref{theo:Ah-approx} applied to $\tebu$ we get that there exists $\tebv$ being $\bopA$-harmonic in $B_{r/2}$ and such that
\[\barint_{B_{r/2}}|\teDbu-\teDbv|\,dx\leq \ve.\]
Then~\eqref{comp-est} follows by rescaling back with $\tev=\lambda\tebv$ which is $\opA$-harmonic.
\end{proof}


\subsection{Estimates on concentric balls} This subsection is devoted to prove some properties of  weak solutions to~\eqref{eq:mu} with $\temu\in C^\infty(\Omega,\Rm)$ holding over a family of concentric balls $\{B^j\}$. Before we pass to this, let us fix some notation and parameters. 
Recall that we have chosen $R_0 = R_0(\data, \varsigma)$ in Proposition~\ref{prop:osc}.
We take an arbitrary constant $\alpha_V \in (0,1)$ and take $\varsigma=\varsigma(s_G,\alpha_V)$ to satisfy $\alpha_{D}:=\frac{\alpha_V+1}{2} \leq 1+(\varsigma-1)/s_G$.
To prove Theorem~\ref{theo:pointwise}, it is enough to take $\alpha_V=\frac{1}{2}$, but for the later use in the proof of Theorem~\ref{theo:H-cont}, we have taken $\alpha_V$ arbitrarily.
We now choose\begin{equation}
    \label{sigma}\sigma_0:= \min\left\{ \left(\frac{1}{2^{n+6} c_{\rm o}}\right)^{\frac{1-\alpha_V}{2}},\frac{1}{4}\right\}.
\end{equation}

If $r\in(0,R_0)$ is given, for every $j\in \N\cup\{0\}$ let us fix
\[r_j:=\sigma^{j+1}r,\qquad\qquad B^j:=\overline{B_{r_j}(x_0)},\]
so that $r_{-1}=r.$ We denote\begin{equation}
    \label{AjEj}
    E_j:=\barint_{B^j}|\teu-(\teu)_{B^j}|\,dx.
\end{equation}


\begin{lem}
 \label{lem:8.1}Suppose Assumption {\bf (A-vect)} is satisfied. If $\teu\in W^{1,G}(\Omega,\Rm)$ is a~weak solution to~\eqref{eq:mu} with $\temu\in C^\infty(\Omega,\Rm)$, $j\in\N$ is fixed, $E_j$ is given by~\eqref{AjEj}, while $0<\sigma\leq \sigma_0$ is arbitrary, then we have that\begin{equation}
     \label{8.12} E_{j+1}\leq c_{\rm D}\sigma^{\alpha_D}E_j+c_{\rm E}r_jg^{-1}\left(\frac{|\temu|(B^j)}{r_j^{n-1}}\right)
 \end{equation} for $c_{\rm D}=c_{\rm D}(\data,\alpha_V) =2^{n+4}c_{\rm o}$ and $c_{\rm E}=c_{\rm E}(\data,\alpha_V)={2^{n+2}c_{\rm s}c_{\rm P}c_{\rm o}}\sigma^{-n},$ where $c_{\rm P}$ is the constant from Poincar\'e inequality in $W^{1,1}(\Omega,\Rm)$.
\end{lem}
\begin{proof} We may apply Lemma~\ref{lem:A-h-appr} in $B_r=B_{r_j}(x_0)$ to get that there exists an $\opA$-harmonic map $\tev_j\in W^{1,G}(B_{r_{j/2}},\Rm)$ in $B_{r_{j/2}}$ such that\begin{equation}
     \label{8.13} \barint_{\frac{1}{2}B^j}|\teDu-\teDv_j|\,dx\leq\frac{\ve}{r_j}{E_j}+c_{\rm s} g^{-1}\left(\frac{|\temu|(B^j)}{r_j^{n-1}}\right).
 \end{equation}
By Poincar\'e inequality in $W^{1,1}(\Omega,\Rm)$ for $\tew_j=\teu-\tev_j$ we have
\[\barint_{\frac{1}{2}B^j}|\tew_j-(\tew_j)_{\frac{1}{2}B^j}|\,dx\leq c_{\rm P} r_j \,\barint_{\frac{1}{2}B^j}|\teDu-\teDv_j|\,dx.\]
Then
\begin{equation}
    \label{est-diff-w_j}
\barint_{\frac{1}{2}B^j}|\tew_j-(\tew_j)_{\frac{1}{2}B^j}|\,dx\leq {\ve}c_{\rm P}\,E_j+c_{\rm s} c_{\rm P}\,r_j g^{-1}\left(\frac{|\temu|(B_{r_j})}{r_j^{n-1}}\right).
\end{equation}
Thus by Lemma~\ref{lem:excess}, the triangle inequality,~\eqref{est-diff-w_j}, and Proposition~\ref{prop:osc},  we estimate  
\begin{flalign*} 
E_{j+1}&= 
\barint_{B^{j+1}}|\teu-(\teu)_{B^{j+1}}|\,dx\\
&\leq \,\barint_{B^{j+1}}|\tew_{j}-(\tew_{j})_{B^{j+1}}|\,dx+\, \barint_{B^{j+1}}|\tev_{j}-(\tev_{j})_{B^{j+1}}|\,dx\\
&\leq \,2\,\barint_{B^{j+1}}|\tew_{j}-(\tew_{j})_{\frac 12 B^{j}}|\,dx+\, 2 c_{\rm o}\sigma^{1+(\varsigma - 1)/s_{G}}\,\barint_{\frac{1}{2}B^{j}}|\tev_{j}-(\tev_{j})_{\frac{1}{2}B^{j}}|\,dx\\
&\leq \left(\frac{2^{n+1}}{\sigma^{n}} + 2c_{\rm o}\sigma^{\alpha_V} \right) \, \barint_{\frac{1}{2}B^{j}}|\tew_j-(\tew_j)_{\frac{1}{2}B^j}|\,dx+\, 2 c_{\rm o} \sigma^{\alpha_V} \barint_{\frac{1}{2}B^{j}}|\teu-(\teu)_{\frac{1}{2}B^{j}}|\,dx\\
&\leq 2^{n+2} c_{\rm o}\sigma^{\alpha_V} E_{j} + \left(\frac{2^{n+1}}{\sigma^{n}} + 2 c_{\rm o}\sigma^{\alpha_V} \right) \,\barint_{\frac{1}{2}B^{j}}|\tew_j-(\tew_j)_{\frac{1}{2}B^j}|\,dx\\
&\leq \left(2^{n+2} c_{\rm o}\sigma^{\alpha_V} + 2 \ve c_{\rm o} c_{\rm P} \sigma^{\alpha_V} + \ve c_{\rm P}\frac{2^{n+1}}{\sigma^{n}} \right)\,E_j\\& \quad +c_{\rm s} c_{\rm P} \left(\frac{2^{n+1}}{\sigma^{n}} + 2 c_0\sigma^{\alpha_V} \right) \,r_j g^{-1}\left(\frac{|\temu|(B_{r_j})}{r_j^{n-1}}\right).
\end{flalign*}
By choosing  $\ve=\frac{\sigma^{n+\alpha_V}}{c_{\rm P}}$ we complete the proof.
\end{proof}

\begin{prop}\label{prop:comp-exc} Suppose Assumption {\bf (A-vect)} is satisfied. If $\teu\in W^{1,G}(\Omega,\Rm)$ is a weak solution to~\eqref{eq:mu} with $\temu\in C^\infty(\Omega,\Rm)$, then there exists a constant $c_V=c_V(\data,\alpha_V)\geq 1$ such that for every $\tau\in(0,1]$ we have\begin{flalign*}
    \barint_{B_{\tau r}(x_0)}&|\teu-(\teu)_{B_{\tau r(x_0)}}|\,dx\\
    &\leq c_V\tau^{\alpha_V}\barint_{B_r(x_0)}|\teu-(\teu)_{B_r(x_0)}|\,dx+c_V\sup_{0<\varrho<r}\vr\, g^{-1}\left(\frac{|\temu|(B_\vr(x_0))}{\vr^{n-1}}\right).
\end{flalign*}
\end{prop}\begin{proof} Lemma~\ref{lem:8.1}  implies that for $j\in\N\cup\{0\}$ it holds that
\[E_{j+1}\leq c_{\rm D}\sigma^{\alpha_D}E_j+c_{\rm E}r_jg^{-1}\left(\frac{|\temu|(B^j)}{r_j^{n-1}}\right).\]
Iterating this estimate we get that for any $k \in\N\cup\{0\}$
\[E_{k+1}\leq (c_{D}\sigma^{\alpha_D})^{k+1} E_0+c_{\rm E} \sum_{j=0}^k (c_{\rm D} \sigma^{\alpha_{\rm D}})^{j}r_j\,g^{-1}\left(\frac{|\temu|(B^j)}{r_j^{n-1}}\right),\]
where $c=c(\data)$.
Recalling $\alpha_{D}=\frac{\alpha_V+1}{2}$, $c_{\rm D}=2^{n+4}c_{\rm o}$ and \eqref{sigma}, we see 
\[ c_{\rm D} \sigma^{\alpha_{\rm D}} \leq \frac{\sigma^{\alpha_V}}{4}. \]
By Lemma~\ref{lem:excess} and direct computation we have for any $k\in\N\cup\{0\}$ that\begin{equation}\label{Ek-est}
    E_k\leq \sigma^{k \alpha_V}\barint_{B_r(x_0)}|\teu-(\teu)_{B_r(x_0)}|\,dx + 2 c_{\rm E}\sup_{0<\vr<r}\vr\,g^{-1}\left(\frac{|\temu|(B_{\vr}(x_0))}{\vr^{n-1}}\right).
\end{equation}
We take $\tau\in(0,\sigma)$ and $k\geq 1$ such that $\sigma^{k+1}<\tau\leq\sigma^k.$ Then by Lemma~\ref{lem:excess} and \eqref{Ek-est} we obtain \begin{flalign*}
& \barint_{B_{\tau r}(x_0)}|\teu-(\teu)_{B_{\tau r}(x_0)}|\,dx \\
&\leq\frac{2\sigma^{kn}}{\tau^n} \barint_{B_{r\sigma^k}(x_0)}|\teu-(\teu)_{B^{k-1}}|\,dx\leq \frac{2E_{k-1}}{\sigma^n}\\
&\leq \frac{\sigma^{(k+1)\alpha_V}}{\sigma^{n+2\alpha_V}} \barint_{B_r(x_0)}|\teu-(\teu)_{B_r(x_0)}|\,dx + \frac{2 c_{\rm E}}{\sigma^n}\sup_{0<\vr<r}\vr\,g^{-1}\left(\frac{|\temu|(B_{\vr}(x_0))}{\vr^{n-1}}\right)\\
&\leq \frac{\tau^{\alpha_V}}{\sigma^{n+2 \alpha_V}} \barint_{B_r(x_0)}|\teu-(\teu)_{B_r(x_0)}|\,dx +\frac{2 c_{\rm E}}{\sigma^n}\sup_{0<\vr<r}\vr\,g^{-1}\left(\frac{|\temu|(B_{\vr}(x_0))}{\vr^{n-1}}\right).
\end{flalign*}
By taking $c_V=c_V(\data,\alpha_V)=2c_{\rm E}\sigma^{-n-2\alpha_V}$ we conclude the claim for $\tau\in(0,\sigma).$ For completing the range of $\tau\in[\sigma,1]$ it suffices to note that
\[\barint_{B_{\tau r}(x_0)}|\teu-(\teu)_{B_{\tau r}(x_0)}|\,dx\leq \frac{2}{\sigma^n}
\barint_{B_{r}(x_0)}|\teu-(\teu)_{B_{r}(x_0)}|\,dx.\]
\end{proof}


\subsection{SOLA $\teu$ belongs to VMO}
\begin{proof}[Proof of Proposition~\ref{prop:vmo}]   Suppose $\temu\in C^\infty(\Omega,\Rm)$ and $\teu\in W^{1,G}(\Omega,\Rm)$ is a weak solution to~\eqref{eq:mu}.
By Proposition~\ref{prop:comp-exc} we can find constants $c_V=c_V(\data,\alpha_V)$ we have\begin{equation}\label{est-weak}
    \barint_{B_{\tau r}}|\teu-(\teu)_{B_{\tau r}}|\,dx\leq c_V\tau^{\alpha_V}\barint_{B_r}|\teu-(\teu)_{B_r}|\,dx+c_V\sup_{0<\varrho<r}\vr\, g^{-1}\left(\frac{|\temu|(B_\vr)}{\vr^{n-1}}\right).
\end{equation}  Let us consider a SOLA $\teu\in W^{1,g}(\Omega,\Rm) $ existing due to Proposition~\ref{prop:exist}. Suppose $(\teu_h)$ and $(\temu_h)$ are approximating sequences from definition of SOLA, see Section~\ref{ssec:main-results}. Inequality~\eqref{est-weak} hold for each  $\teu_h$ and $\temu_h$. We have to motivate passing to the limit with $h\to\infty$. Since~\eqref{conv-of-meas} holds, we can write~\eqref{est-weak}  for the original SOLA too. From now on this kind of solution is considered.

 Our aim now is to show that $\teu$ is VMO at $x_0$ provided \eqref{mu-shrinks} is assumed. Let $\delta\in(0,1)$. By~\eqref{mu-shrinks} we find a positive radius $r_{1,\delta}<r$ such that \[c_V\sup_{0<\vr<r_{1,\delta}}\vr\, g^{-1}\left(\frac{|\temu|(B_\vr(x_0))}{\vr^{n-1}}\right)\leq \frac{\delta}{2}\]
and then $\tau_\delta$ so small that 
\[c_V\tau_\delta^{\alpha_V}\barint_{B_{r_{1,\delta}}}|\teu-(\teu)_{B_{r_{1,\delta}}(x_0)}|\,dx\leq\frac{\delta}{2}.\]
For $r_\delta:=\tau_\delta r_{1,\delta}$ from estimate~\eqref{est-weak} (applied with $r=r_{1,\delta}$) it follows that 
\[\sup_{0<\vr<r_\delta}\barint_{B_{\vr}(x_0)}|\teu-(\teu)_{B_{\vr(x_0)}}|dx\leq\delta,\]
that is that $\teu$ has vanishing mean oscillation at $x_0$.
\end{proof}


\subsection{Proofs of Theorems~\ref{theo:pointwise},~\ref{theo:continuity} and~\ref{theo:H-cont}} We start with the proof of pointwise Wolff potential estimate, then pass to continuity and H\"older continuity criteria.
\begin{proof}[Proof of Theorem~\ref{theo:pointwise}]
We notice that having $E_j$ defined in~\eqref{AjEj} with $r=r^j$ we can fix $\sigma$ in Lemma~\ref{lem:8.1} to get that\begin{equation}
    \label{est-2-st-1-theo-est}
    E_{j+1}\leq \frac{1}{2}E_j+cr^jg^{-1}\left(\frac{|\temu|(B^j)}{r_j^{n-1}}\right)\qquad\text{for every }\ j\in\N\cup\{0\}.
\end{equation}
We sum up inequalities from~\eqref{est-2-st-1-theo-est} to obtain
\[\sum_{j=1}^{k+1} E_j\leq\frac{1}{2}\sum_{j=0}^k E_j+c\sum_{j=0}^k r_jg^{-1}\left(\frac{|\temu|(B^j)}{r_j^{n-1}}\right),\qquad k\in\N\cup\{0\}.\]
By rearranging terms we have
\[\sum_{j=1}^{k+1} E_j\leq 2 E_0+c\sum_{j=0}^k r_jg^{-1}\left(\frac{|\temu|(B^j)}{r_j^{n-1}}\right).\]
We notice that for some $c=c(\data)$ we can estimate
\[\sum_{j=0}^k r_jg^{-1}\left(\frac{|\temu|(B^j)}{r_j^{n-1}}\right)\leq c\,\int_0^r g^{-1}\left(\frac{|\temu|(B_\vr)}{\vr^{n-1}}\right)\,d\vr=c\,\cW^\temu_G(x_0,r).\]
Last two displays imply that
\[\sum_{j=1}^{k+1} E_j\leq 2 E_0+c\cW^\temu_G(x_0,r).\]
For every $m,k\in\N$ such that $m<k$ we have\begin{flalign*}
|(\teu)_{B^k}-(\teu)_{B^m}|&\leq \sum_{j=m}^{k-1} |(\teu)_{B^{j+1}}-(\teu)_{B^j}|\leq \sigma^{-n} \sum_{j=m}^{k+1}E_j
\leq \sigma^{-n} \sum_{j=0}^{k+1}E_j\\
&\leq 2\sigma^{-n} E_0+c\sigma^{-n}\cW^\temu_G(x_0,r)\\
&\leq 2\sigma^{-n}\barint_{B_r(x_0)}|\teu-(\teu)_{B_r(x_0)}|\,dx+c\sigma^{-n}\cW^\temu_G(x_0,r),
\end{flalign*}
where $\sigma=\sigma(\data)$ and $c=c(\data).$ For $j\to\infty$, $((\teu)_{B^j})_j$ is a Cauchy sequence that converges to $\teu(x_0)$, that is
\[\lim_{\vr\to 0}(\teu)_{B_\vr(x_0)}=\teu(x_0)\]
and $x_0$ is a Lebesgue's point of $\teu$. This completes the proof of \eqref{Wolff-osc-est}, while~\eqref{eq:u-est} follows as a direct corollary.\end{proof}


Let us concentrate on the continuity criterion.
\begin{proof}[Proof of Theorem~\ref{theo:continuity}] 
Our aim is to infer continuity of $\teu$ in $B_r(x_0)$ knowing that~\eqref{Wolff-shrinks} holds. We will show that for every $\delta>0$ and $x_1\in B_r(x_0)$ we can find $r_\delta\in (0,{\rm dist}\,(B_r(x_0),\partial\Omega))$ such that\begin{equation}
    \label{osc-u}{\rm osc}_{B_{r_\delta}(x_1)}\,\teu<\delta.
\end{equation} Without loss of generality we assume that $\temu$ is defined on whole $\rn$, as we can extend it by zero outside $\Omega$. 
By~\eqref{Wolff-shrinks} we can take $\vr_1$ small enough for
\begin{equation}
    \label{small-wolff} \sup_{x\in B_r(x_0)}\cW^\temu_G(x,\vr_1)\leq \frac{\delta}{16}.
\end{equation}
Let $r_\delta>0$ to be chosen in a moment. We take an arbitrary point $x_2\in B_{r_\delta}(x_1)$ and estimate \begin{flalign}\nonumber |\teu(x_1)-\teu(x_2)| &\leq |\teu(x_1)-(\teu)_{B_{2r_\delta}(x_1)}|+|(\teu)_{B_{2r_\delta}(x_1)}-(\teu)_{B_{r_\delta}(x_2)}|\\
&\quad+|(\teu)_{B_{r_\delta}(x_2)}-\teu(x_2)|=:A_1+A_2+A_3.\label{telescope}
\end{flalign}
 We start with estimating $A_2$ by noting that \[A_2=|(\teu)_{B_{r_\delta}(x_2)}-(\teu)_{B_{2r_\delta}(x_1)}| \leq  \barint_{B_{r_\delta}(x_2)}|\teu-(\teu)_{B_{2r_\delta}(x_1)}|\,dx.\] Since~\eqref{Wolff-shrinks} implies~\eqref{mu-shrinks}, Proposition~\ref{prop:vmo} implies that $\teu$ has vanishing mean oscillations at~$x_1$. Therefore there exists $\vr_2\in(0,\min\{\vr_1,{\rm dist}(x_1,\partial B_r(x_0))/4\})$ such that for every $\vr\leq \vr_2$ it holds\[\barint_{B_{\vr}(x_1)}|\teu-(\teu)_{B_{\vr}(x_1)}|\,dx\leq\frac{\delta}{2^{n+4}}.\]
We choose $r_\delta=\vr_2/2$  and observe that~\eqref{small-wolff} imply that\begin{flalign*}\barint_{B_{r_\delta}(x_2)}|\teu-(\teu)_{B_{2 r_\delta}(x_1)}|\,dx
&\leq 2^{n}\,\barint_{B_{2r_\delta}(x_1)}|\teu-(\teu)_{B_{2 r_\delta}(x_1)}|\,dx\leq\frac{\delta}{16}.
\end{flalign*}
In turn $ A_2 \leq\frac{\delta}{16}.$ By Theorem~\ref{theo:pointwise} and~\eqref{small-wolff} we get that $x_1$ and $x_2$ are Lebesgue's points and \[A_1+A_3=|\teu(x_1)-(\teu)_{B_{2r_\delta}(x_1)}|+|(\teu)_{B_{r_\delta}(x_2)}-\teu(x_2)|\leq \frac{\delta}4 .\]
Applying these observation we get from~\eqref{telescope} that\begin{flalign*} |\teu(x_1)-\teu(x_2)|&\leq  \frac{\delta}{2}.
\end{flalign*}
Since $x_2$ was an arbitrary point of $B_{r_\delta}(x_1)$, we have \eqref{osc-u} justified, which completes the proof.
\end{proof}

We are in the position to prove the H\"older continuity criterion.
\begin{proof}[Proof of Theorem~\ref{theo:H-cont}.] 
Notice that assumption~\eqref{mu-control} implies that there exists $c=c(\data)>0$, such that for all sufficiently small $r$ we have\begin{equation*}
    \cW^\temu_G(x,r)\leq c r^\theta.
\end{equation*}
Applying assumption~\eqref{mu-control} to Proposition~\ref{prop:comp-exc} with $\alpha_V=\frac{\theta+1}{2}$, we have
\begin{flalign*}
    \barint_{B_{\rho}} |\teu-(\teu)_{B_{\rho}}|\,dx\leq c \left( \frac{\rho}{r} \right)^{\frac{\theta+1}{2}} \barint_{B_r}|\teu-(\teu)_{B_r}|\,dx + c r^{\theta}
\end{flalign*}
for any $0 < \rho < r \leq R_0$.
Now we apply Lemma~\ref{lem:absorb2} to see
\[ \barint_{B_{\rho}} |\teu-(\teu)_{B_{\rho}}|\,dx\leq c \left( \frac{\rho}{r} \right)^{\theta} \barint_{B_r}|\teu-(\teu)_{B_r}|\,dx + c \rho^{\theta} \]
By Campanato's characterization \cite[Theorem~2.9]{giusti}, we complete the proof.
\end{proof}

\section*{Appendix}
\begin{proof}[Proof of Lemma~\ref{lem:Wolff-est}]Set $x\in\Omega_0,$ $R_k=2^{1-k}R$ and $B_k=B_{R_k}(x)$ for $k=0,1,\dots\,$. As $\teF$ is taken in a place of a measure with a slight abuse of notation we write $|\teF|(B_{R_k}(x))=\int_{B_k}|\teF(y)|\,dy.$ We notice that 
we have
\[\cW_G^{|\teF|}(x,R)=\sum_{k=1}^\infty \int_{R_{k+1}}^{R_k} g^{-1}\left(\frac{|\teF|(B_r(x))}{r^{n-1}}\right)\,dr\lesssim \sum_{k=1}^\infty  R_{k} g^{-1}\left(\frac{|\teF|(B_{R_k}(x))}{R_k^{n-1}}\right)\,.\]
To estimate the series we employ the decreasing rearrangement $|\teF|^\star$ of $|\teF|$ and its maximal rearrangement $|\teF|^{\star\star}$. When $w_n$ is the volume of the unit ball, we have that
\begin{flalign*}
  \frac{|\teF|(B_{R_k}(x))}{R_k^{n-1}}&=\frac{1}{R_k^{n-1}}\int_{B_{R_k}(x)}|\teF(y)|\,dy\\&\leq {w_n R_k}\, \barint_0^{w_nR_k^n}|\teF|^\star(t)\,dt={w_n R_k}\, |\teF|^{\star\star}({w_nR_k^n})\,.
\end{flalign*}
Then  we have \begin{flalign}\label{1na8} R_kg^{-1}\left(\frac{|\teF|(B_{R_k}(x))}{R_k^{n-1}}\right)&\lesssim R_kg^{-1}\left({w_n R_k}\, |\teF|^{\star\star}({w_nR_k^n})\right)\\ &\lesssim  \int_{w_n R_k^{n}}^{w_n R_{k-1}^n}\rho^\frac{1}{n}g^{-1}\left(\rho^\frac{1}{n} \, |\teF|^{\star\star}(\rho)\right)\,\frac{d\rho}{\rho}\nonumber
\end{flalign}
with implicit constants independent of $k$. Therefore
\begin{flalign*}\sup_{x\in\Omega_0}\cW_G^{|\teF|}(x,R)&\lesssim \sum_{k=1}^\infty \int_{w_n R_k^{n}}^{w_n R_{k-1}^n}  \rho^{\frac{1}{n}-1} g^{-1}\left(\rho^\frac{1}{n}\, |\teF|^{\star\star}(\rho)\right)\,\frac{d\rho}{\rho}\\
&= \int_{0}^{w_n R^n} \rho^\frac{1}{n} g^{-1}\left(\rho^\frac{1}{n}\, |\teF|^{\star\star}(\rho)\right)\,\frac{d\rho}{\rho}\,.
\end{flalign*}
\end{proof} 

\noindent{\bf Acknowledgements} I. Chlebicka is supported by NCN grant no. 2019/34/E/ST1/00120. A. Zatorska-Goldstein is supported by NCN grant no. 2019/33/B/ST1/00535. Y. Youn is supported by NRF grant no. 2020R1C1C1A01009760.

\end{document}